\documentclass[aos,MSNbibl,seceqn,dvips]{arximspdf}
\usepackage{algorithm}
\usepackage{dcolumn}
\usepackage{graphicx}

\doi{10.1214/12-AOS1057} 
\volume{40}
\issue{6}
\pubyear{2012}
\firstpage{3077}
\lastpage{3107}

\makeatletter

\newcolumntype{d}[1]{D{.}{.}{#1}}

\newcommand{\rright}{\right}
\newcommand{\lleft}{\left}

\newtheorem{lemma}{Lemma}[section]
\newtheorem{proposition}{Proposition}[section]
\newtheorem{theorem}{Theorem}[section]

\newtheorem{observation}{Observation}[section]

\newcommand{\cal}{\mathcal}

\def\Col{\mathrm{Col}}
\def\Row{\mathrm{Row}}
\def\BRIP{\operatorname{BRIP}}
\def\RE{\operatorname{RE}}

\def\Tr{\operatorname{Tr}}
\def\sign{\operatorname{sign}}

\newcommand{\bbr}{{\mathbb{R}}}
\newcommand{\bbs}{{\mathbb{S}}}
\newcommand{\cW}{{\cal W}}

\def\Ker{\operatorname{Ker}}
\def\Prob{\operatorname{Prob}}

\def\Argmin{\mathop{\operatorname{Arg}\operatorname{min}}}
\def\pen{{\mathrm{pen}}}
\def\reg{{\mathrm{reg}}}
\def\bQ{{\mathbf{Q}}}
\def\bH{{\mathbf{H}}}
\def\cU{{\cal U}}
\def\cS{{\cal S}}
\def\cH{{\cal H}}
\def\bT{{\mathbf{T}}}
\def\ErfInv{\operatorname{Erfinv}}
\def\sErfInv{\operatorname{Erfinv}}
\def\cN{{\cal N}}
\def\bS{{\mathbf{S}}}
\def\bR{{\mathbf{R}}}

\makeatother

\begin{document}
\begin{frontmatter}

\title{Accuracy guaranties for $\ell_1$ recovery of
block-sparse signals}
\runtitle{Accuracy guaranties for block-$\ell_1$ recovery}

\begin{aug}
\author[A]{\fnms{Anatoli} \snm{Juditsky}\corref{}\ead[label=e1]{juditsky@imag.com}},
\author[B]{\fnms{Fatma}~\snm{K{\i}l{\i}n\c{c} Karzan}\thanksref{t1}\ead[label=e2]{fkilinc@andrew.cmu.edu}},
\author[C]{\fnms{Arkadi}~\snm{Nemirovski}\thanksref{t1,t2}\ead[label=e3]{nemirovs@isye.gatech.edu}}
\and
\author[D]{\fnms{Boris} \snm{Polyak}\ead[label=e4]{boris@ipu.rssi.ru}}
\runauthor{Juditsky, K{\i}l{\i}n\c{c} Karzan, Nemirovski and Polyak}
\affiliation{Universit\'e J. Fourier de Grenoble, Carnegie Mellon
University,
Georgia~Institute of Technology and Institute of Control Sciences}
\address[A]{A. Juditsky\\
LJK\\
Universit\'e J. Fourier\\
B.P. 53\\
38041 Grenoble Cedex 9\\
France\\
\printead{e1}}
\address[B]{F. K{\i}l{\i}n\c{c} Karzan\\
Tepper School of Business\\
Carnegie Mellon University\\
Pittsburgh, Pennsylvania 15213\hspace*{9pt}\\
USA\\
\printead{e2}}
\address[C]{A. Nemirovski\\
School of Industrial and Systems Engineering \\
Georgia Institute of Technology\\
Atlanta, Georgia 30332\\
USA\\
\printead{e3}}
\address[D]{D. Polyak\\
Institute of Control Sciences\\
\quad of Russian Academy of Sciences\\
Moscow 117997\\
Russia\\
\printead{e4}} 
\end{aug}

\thankstext{t1}{Supported by the Office of Naval Research Grant N000140811104.}
\thankstext{t2}{Supported by the NSF Grant DMS-09-14785.}

\received{\smonth{11} \syear{2011}}
\revised{\smonth{9} \syear{2012}}

%
\begin{abstract}
We introduce a general framework to handle structured models (sparse
and block-sparse with possibly overlapping blocks). We discuss new
methods for their recovery from incomplete observation, corrupted with
deterministic and stochastic noise, using block-$\ell_1$
regularization. While the current theory provides promising bounds for
the recovery errors under a number of different, yet mostly hard to
verify conditions, our emphasis is on verifiable conditions on the
problem parameters (sensing matrix and the block structure) which
guarantee accurate recovery. Verifiability of our conditions not only
leads to efficiently computable bounds for the recovery error but also
allows us to optimize these error bounds with respect to the method
parameters, and therefore construct estimators with improved
statistical properties. To justify our approach, we also provide an
oracle inequality, which links the properties of the proposed recovery
algorithms and the best estimation performance. Furthermore, utilizing
these verifiable conditions, we develop a computationally cheap
alternative to block-$\ell_1$ minimization, the non-Euclidean Block
Matching Pursuit algorithm. We close by presenting a numerical study to
investigate the effect of different block regularizations and
demonstrate the performance of the proposed recoveries.
\end{abstract}

%

\end{frontmatter}

\section{Introduction}\label{secintro}

\subsection*{The problem} Our goal in this paper is to estimate a
linear transform $Bx \in\bbr^N$ of a vector $x\in\bbr^n$ from the
observations
%
\begin{equation}
\label{anobs} y=Ax+u+\xi.
\end{equation}
Here $A$ is a given $m\times n$ sensing matrix, $B$ is a given $N\times
n$ matrix, and $u+\xi$ is the observation error; in this error, $u$ is
an unknown \textit{nuisance} known to belong to a given compact convex set
$\cU\subset\bbr^m$ symmetric w.r.t. the origin, and $\xi$ is random
noise with known distribution $P$.

We assume that the space $\bbr^N$ where $Bx$ lives is represented as
$\bbr^N=\bbr^{n_1}\times\cdots\times\bbr^{n_K}$, so that a vector
$w\in\bbr^N$ is a block vector: $w=[w[1];\ldots;w[K]]$ with blocks
$w[k]\in\bbr^{n_k}$, $1\leq k\le K$.\setcounter{footnote}{2}\footnote{We use MATLAB
notation: $[u,v,\ldots,z]$ is the horizontal concatenation of matrices
$u,v,\ldots,z$ of common height, while $[u;v;\ldots;z]$ is the vertical
concatenation of matrices $u,v,\ldots,z$ of common width. All vectors are
column vectors.} In particular, $Bx=[B[1]x;\ldots;B[K]x]$ with $n_k\times
n$ matrices $B[k]$, $1\leq k\leq K$. While we do not assume that the
vector $x$ is sparse in the usual sense, we do assume that the linear
transform $Bx$ to be estimated is \textit{$s$-block sparse}, meaning that at
most a given number, $s$, of the blocks $B[k]x$, $1\leq k\leq K$, are nonzero.

The recovery routines we intend to consider are based on \textit{block-$\ell_1$ minimization}, that is, the estimate $\widehat{w}(y)$
of $w=Bx$ is $B\widehat{z}(y)$, where $\widehat{z}(y)$ is obtained by
minimizing the norm ${\sum_{k=1}^K}\|B[k]z\|_{(k)}$ over signals $z\in
\bbr^n$ with $Az$ ``fitting,'' in a certain precise sense, the
observations $y$. Above, $\|\cdot\|_{(k)}$ are given in advance norms
on the spaces $\bbr^{n_k}$ where the blocks of $Bx$ take their values.

In the sequel we refer to the given in advance collection $\cS
=(B,n_1,\ldots,n_K,\break {\|\cdot\|_{(1)}},\ldots,\|\cdot\|_{(K)})$ as the
\textit{representation structure} (r.s.). Given such a representation structure
$\cS$ and a sensing matrix $A$, our ultimate goal is to understand how
well one can recover the $s$-block-sparse transform $Bx$ by
appropriately implementing block-$\ell_1$ minimization.

\subsection*{Related Compressed Sensing research} Our situation and
goal form a straightforward extension of the usual sparsity/block
sparsity framework of Compressed Sensing. Indeed,
the \textit{standard representation structure} with $B=I_n$, $n_k=1$, and
$\|\cdot\|_{(k)}=|\cdot|$, $1\leq k\leq K=n$, leads to the standard
Compressed Sensing setting---recovering a sparse signal $x\in\bbr^n$
from its noisy observations (\ref{anobs}) via $\ell_1$ minimization.
The case of nontrivial block structure $\{n_k,\|\cdot\|_{(k)}\}
_{k=1}^K$ and $B=I_n$ is generally referred to as \textit{block-sparse},
and has been considered in numerous recent papers.
Block-sparsity (with $B=I_n$) arises naturally (see, e.g., \cite
{EldarMishali09} and references therein) in a number of applications
such as multi-band signals, measurements of gene expression levels or
estimation of multiple measurement vectors sharing a joint sparsity pattern.
Several methods of estimation and selection extending the ``plain''
$\ell_1$-minimization to block sparsity were proposed and investigated
recently. Most of the related research focused so far on \textit{block
regularization schemes}---group Lasso recovery of the form
\[
\widehat{x}(y)\in\Argmin_{z=[z^1;\ldots;z^K]\in\bbr^n=\bbr
^{n_1}\times\cdots\times\bbr^{n_K}} \Biggl\{\|Az-y\|_2^2+
\lambda\sum_{k=1}^K\bigl\|z[k]\bigr\|_2
\Biggr\}
\]
(here \mbox{$\|\cdot\|_2$} is the Euclidean norm of the block).
In particular, the
literature on ``plain Lasso'' (the case of $n_k=1, 1\leq k\leq K=n$)
%
%
has an important counterpart on group Lasso; see, for example, \cite
{Bach08GroupLasso,Ben-HaimEldar10,ChesneauHebiri08,DuarteBajwaCalderbank11,EldarKuppingerBolcskei10,EldarMishali09,GribonvalNielsen03,HuangZhang10,LiuZhang09,MeiervandeGeerBuhlmann08,NardiRinaldo08,Obozinskietal11,VikaloParvaresh07,StojnicParvareshHassibi09,YuanLin06}
and the references therein. Another celebrated technique of sparse
recovery, the Dantzig selector, originating from \cite{CandesTao07},
has also been extended to handle block-sparse structures \cite
{JamesRadchenkoLv09,groupDantzig10}. Most of the cited papers focus on
bounding recovery errors in terms of the magnitude of the observation
noise and ``$s$-concentration'' of the true signal~$x$ (the distance
from the space of signals with at most $s$ nonzero blocks---the sum of
magnitudes $\|x[k]\|_2$ of all but the $s$ largest in magnitude blocks
in $x$). Typically, these results rely on natural block analogy
(``Block RIP;'' see, e.g., \cite{EldarMishali09}) of the celebrated
Restricted Isometry Property introduced by Cand\'{e}s and Tao
\mbox{\cite{CandesTaorip05,Candes08note}} or on block analogies
\cite{Lounicietal10} of the Restricted Eigenvalue Property introduced
in~\cite{BickelRitovTsybakov08}. In addition to the usual
(block)-sparse recovery, our framework also allows to handle group
sparse recovery with overlapping groups by properly defining the
corresponding $B$ matrix.

\subsection*{Contributions of this paper} The first (by itself, minor)
novelty in our problem setting is the presence of the linear mapping
$B$. We are not aware of any preceding work handling the case of a
``nontrivial'' (i.e., different from the identity) $B$. We qualify this
novelty as minor, since in fact the case of a nontrivial $B$ can be
reduced to the one of $B=I_n$.\footnote{Assuming, for example, that
$x\mapsto Bx$ is an ``onto'' mapping, we can treat $Bx$ as our signal,
the observations being $Py$, where $P$ is the projector onto the
orthogonal complement to the linear subspace $A\cdot\Ker B$ in $\bbr
^m$; with $y=Ax+u+\xi$, we have $Py=GBx+P(u+\xi)$ with an explicitly
given matrix $G$.} However, ``can be reduced'' is not the same as
``should be reduced,'' since problems with nontrivial $B$ mappings
arise in many applications. This is the case, for example, when $x$ is
the solution of a linear finite-difference equation with a sparse
right-hand side (``evolution of a linear plant corrected from time to
time by impulse control''), where $B$ is the matrix of the
corresponding finite-difference operator. Therefore, introducing $B$
adds some useful flexibility (and as a matter of fact costs nothing, as
far as the theoretical analysis is concerned).

We believe, however, that the major novelty in what follows is the
emphasis on \textit{verifiable} conditions on matrix $A$ and the r.s. $\cS
$ which \textit{guarantee} good recovery of the transform $Bx$ from noisy
observations of $Ax$, provided that the transform in question is nearly\vadjust{\goodbreak}
$s$-block sparse, and the observation noise is low. Note that such
efficiently verifiable guarantees cannot be obtained from the
``classical'' conditions\footnote{Note that it has been recently
proved in \cite{PfetschTillmann12} that computing the parameters
involved in verification of Nullspace condition as well as RIP for
sparse recovery is NP-hard.} used when studying theoretical properties
of block-sparse recovery (with a notable exception of the Mutual
Block-Incoherence condition of \cite{EldarKuppingerBolcskei10}).
For example, given $A$ and $\cS$, one cannot answer in any reasonable
time if the (Block-) Restricted Isometry or Restricted Eigenvalue
property holds with given parameters. While the efficient verifiability
is by no means necessary for a condition to
be meaningful and useful, we believe that verifiability has its value
and is worthy of being investigated. In particular,
it allows us to design new recovery routines with explicit confidence
bounds for the recovery error and then optimize these bounds with
respect to the method parameters.
In this respect, the current work extends the results of \cite
{JNCS,JKNCS,JNnoisy}, where $\ell_1$ recovery of the ``usual'' sparse
vectors was considered (in the first two papers---in the case of
uncertain-but-bounded observation errors, and in the third---in the
case of Gaussian observation noise). Specifically, we propose here new
routines of block-sparse recovery which explicitly utilize a \textit{contrast matrix}, a kind of ``validity certificate,'' and show how
these routines may be tuned to attain the best performance bounds. In
addition to this, verifiable conditions pave the way of efficiently
designing sensing matrices which possess certifiably good recovery
properties for block-sparse recovery (see \cite{JKNCSDesign} for
implementation of such an approach in the usual sparsity setting).

The main body of the paper is organized as follows: in Section \ref
{secprob} we formulate the block-sparse recovery problem and introduce
our core assumption---a family of conditions $\bQ_{s,q}$, $1\le q\le
\infty$, which links the representation structure $\cal S$ and sensing
matrix $A\in\bbr^{m\times n}$ with a contrast matrix $H\in\bbr^{m\times M}$.
Specifically, given $s$ and $q\in[1,\infty]$ and a norm $\|\cdot\|$,
the condition $\bQ_{s,q}$ on an $m\times M$ \textit{contrast matrix} $H$
requires $\exists\kappa\in[0,1/2)$ such that
\[
\forall\bigl(x\in\bbr^n\bigr)\qquad  L_{s,q}(Bx) \leq
s^{{1/q}}\bigl\|H^TAx\bigr\| +\kappa s^{{1/ q}-1}
L_1(Bx)
\]
holds, where for
$w=[w[1];\ldots;w[K]]\in\bbr^N$ and $p\in[1,\infty]$,
\[
L_p(w)=\bigl\|\bigl[\bigl\|w[1]\bigr\|_{(1)};\ldots;\bigl\|w[K]\bigr\|_{(K)}\bigr]\bigr\|_p
\]
and
\[
L_{s,p}(w)=\bigl\|
\bigl[\bigl\|w[1]\bigr\|_{(1)};\ldots;\bigl\|w[K]\bigr\|_{(K)}\bigr]\bigr\|_{s,p},
\]
where $\|u\|_{s,p}$ is
the norm on $\bbr^K$ defined as follows: we zero out all but the $s$
largest in magnitude entries in vector $u$, and take the $\|\cdot\|
_p$-norm of the resulting $s$-sparse vector.
Then, by restricting our attention to the standard representation
structures, we study the relation between condition $\bQ_{s,q}$ and
the usual assumptions used to validate block-sparse recovery, for
example, Restricted Isometry/Eigenvalue Properties and their block versions.

In Section~\ref{secrecovery} we introduce two recovery routines based
on the $L_1(\cdot)$ norm:
\begin{itemize}
\item \textit{regular $\ell_1$ recovery} [cf. (block-) Dantzig selector]
\[
\widehat{x}_\reg(y) \in\Argmin_{z\in\bbr^n} \bigl\{
L_1(Bz)\dvtx  \bigl\|H^T(y-Az)\bigr\|_\infty\leq\rho\bigr\},
\]
where with probability $1-\varepsilon$, $\rho[\mbox{$=$}\rho(H,\varepsilon)]$ is
an upper bound on the $\|\cdot\|$-norm of the observation error;
\item \textit{penalized $\ell_1$ recovery} [\textit{cf.} (\textit{block}-) \textit{Lasso}]
\[
\widehat{x}_\pen(y)\in\Argmin_{z\in\bbr^n} \bigl[L_1(Bz)
+ 2s\bigl\|H^T(y-Az)\bigr\|_\infty\bigr],
\]
where $s$ is our guess for the number of nonvanishing blocks in the
true signal~$Bx$.
\end{itemize}
Under condition $\bQ_{s,q}$, we establish performance guarantees of
these recoveries, that is, explicit upper bounds on the size of
confidence sets for the recovery error $L_p(B(\widehat{x}-x))$, $1\le
p\le q$. Our performance guarantees have the usual natural
interpretation---as far as recovery of transforms $Bx$ with small
$s$-block concentration\footnote{$s$-block concentration of a block
vector $w$ is defined as $L_1(w)-L_{s,1}(w)$.} is concerned, everything
is as if we were given the direct observations of $Bx$ contaminated by
noise of small $L_\infty$ magnitude.

Similar to the usual assumptions from the literature, conditions $\bQ
_{s,q}$ are generally computationally intractable, nonetheless, we
point out a notable exception in Section \ref
{secconditionproperties}. When all block norms are $\|\cdot\|_{(k)}=\|
\cdot\|_\infty$, the condition $\bQ_{s,\infty}$, the
strongest among our family of conditions, is efficiently verifiable.
Besides, in this situation, the latter condition is ``fully
computationally tractable,'' meaning that one can optimize efficiently
the bounds for the recovery error over the contrast matrices $H$
satisfying $\bQ_{s,\infty}$ to design optimal recovery routines. In
addition to this, in Section~\ref{secconditionnecessity}, we
establish an oracle inequality which shows that existence of the
contrast matrix $H$ satisfying condition $\bQ_{s,\infty}$ is not only
sufficient but also necessary for ``good recovery'' of block-sparse
signals in the $L_\infty$-norm when $\|\cdot\|_{(k)}=\|\cdot\|_\infty$.

In Section~\ref{secconditionsufficiency} we provide a verifiable
sufficient condition for the validity of $\bQ_{s,q}$ for general $q$,
assuming that $\cS$ is $\ell_r$-r.s. [i.e., $\|\cdot\|_{(k)}=\|\cdot
\|_r$, $1\leq k\leq K$], and, in addition, $r\in\{1,2,\infty\}$. This
sufficient condition can be used to build a ``quasi-optimal'' contrast
matrix $H$. We also relate this condition to the Mutual
Block-Incoherence condition of \cite{EldarKuppingerBolcskei10}
developed for the case of $\ell_2$-r.s. with $B=I_n$. In particular,
we show in Section~\ref{secmutual} that the Mutual Block-Incoherence
is more conservative than our verifiable condition, and thus is
``covered'' by the latter. ``Limits of performance'' of our verifiable
sufficient conditions are investigated in Section~\ref
{seclimitsofperformance}.

In Section~\ref{secNEMP} we describe a computationally cheap
alternative to block-$\ell_1$ recoveries---a non-Euclidean Block
Matching Pursuit (NEBMP) algorithm. Assuming that $\cS$ is either
$\ell_2$-, or $\ell_\infty$-r.s. and that the verifiable sufficient
condition $\bQ_{s,\infty}$ is satisfied, we show that this algorithm
(which does not require optimization) provides performance guarantees
similar to those of regular/penalized $\ell_1$ recoveries.

We close by presenting a small simulation study in Section~\ref{secnumerics}.

Proofs of all results are given in the supplementary
article \cite{JKNP-suppl}.

\section{Problem statement}\label{secprob}

\subsection*{Notation} In the sequel, we deal with:
\begin{itemize}
\item \textit{signals}---vectors $x=[x_1;\ldots;x_n]\in\bbr^n$, and an
$m\times n$ \textit{sensing matrix} $A$;
\item \textit{representations of signals}---block vectors
$w=[w[1];\ldots;w[K]]\in\cW:=\break\bbr^{n_1}_{w[1]}\times\cdots\times\bbr
^{n_K}_{w[K]}$, and the \textit{representation matrix}
$B=[B[1];\ldots;B[K]]$, $B[k]\in\bbr^{n_k\times n}$; the representation
of a signal $x\in\bbr^n$ is the block vector $w=Bx$ with the blocks
$B[1]x,\ldots,B[K]x$.
\end{itemize}
From now on, the dimension of $\cW$ is denoted by $N$:
\[
N=n_1+\cdots+n_K.
\]
The factors $\bbr^{n_k}$ of the representation space $\cW$ are
equipped with norms $\|\cdot\|_{(k)}$; the conjugate norms are denoted
by \mbox{$\|\cdot\|_{(k,*)}$}. A vector $w=[w[1];\ldots;\break w[K]]$ from $\cW
$ is called \textit{$s$-block-sparse}, if the number of nonzero blocks
$w[k]\in\bbr^{n_k}$ in $w$ is at most $s$. A vector $x\in\bbr^n$
will be called $s$-block-sparse, if its representation $Bx$ is so.
We refer to the collection $\cS=(B,n_1,\ldots,n_K,\|\cdot\|_{(1)},\ldots,\|
\cdot\|_{(K)})$ as the \textit{representation structure} (r.s. for
short). The \textit{standard} r.s. is given by $B=I_n$, $K=N$,
$n_1=\cdots=n_n=1$ and $\|\cdot\|_{(k)}=|\cdot|$, $1\leq k\leq N$, and
\textit{an $\ell_r$-r.s.} is the r.s. with $\|\cdot\|_{(k)}=\|\cdot\|_r$,
$1\leq k\leq K$.

For $w\in\cW$, we call the number $\|w[k]\|_{(k)}$ the \textit{magnitude} of the $k$th block in $w$ and denote by $w^s$ the
representation vector obtained from $w$ by zeroing out all but the $s$
largest in magnitude blocks in $w$ (with the ties resolved
arbitrarily). For $I\subset\{1,\ldots,K\}$ and a representation vector
$w$, $w_I$ denotes the vector obtained from $w$ by keeping intact the
blocks $w[k]$ with $k\in I$ and zeroing out all remaining blocks. For
$w\in\cW$ and $1\leq p\leq\infty$, we denote by $L_p(w)$ the $\|
\cdot\|_p$-norm of the vector $[\|w[1]\|_{(1)};\ldots;\|w[K]\|_{(K)}]$,
so that $L_p(\cdot)$ is a norm on $\cW$ with the conjugate norm
$L_p^*(w)=\|[\|w[1]\|_{(1,*)};\ldots;\|w[K]\|_{(K,*)}]\|_{p_*}$
where $p_*=\frac{p}{p-1}$. Given a positive integer $s\leq K$, we set
$L_{s,p}(w)=L_p(w^s)$. Note that $L_{s,p}(\cdot)$ is a norm on $\cW$.
We define the \textit{$s$-block concentration} of a vector $w$ as
$\upsilon_s(w)=L_1(w-w^s)$.

\subsection*{Problem of interest} Given an observation
%
\begin{equation}
\label{newobs} y=Ax+u+\xi,\vadjust{\goodbreak}
\end{equation}
of unknown signal $x\in\bbr^n$, we want to recover the representation
$Bx$ of $x$, knowing in advance that this representation is ``nearly
$s$-block-sparse,'' that is, the representation can be approximated by
an $s$-block-sparse one; the $L_1$-error of this approximation, that
is, the $s$-block concentration, $\upsilon_s(Bx)$, will be present in
our error bounds.

In (\ref{newobs}) the term $u+\xi$ is the observation error; in this
error, $u$ is an unknown \textit{nuisance} known to belong to a given
compact convex set $\cU\subset\bbr^m$ symmetric w.r.t. the origin,
and $\xi$ is random noise with known distribution~$P$.

\subsection*{Condition $\bQ_{s,q}(\kappa)$}

We start with introducing the condition which will be instrumental in
all subsequent constructions and results.
Let a sensing matrix $A$ and an r.s. $\cS=(B,n_1,\ldots,n_K,\|\cdot\|
_{(1)},\ldots,\|\cdot\|_{(K)})$ be given, and let $s\leq K$ be a positive
integer, $q\in[1,\infty]$ and $\kappa\geq0$. We say that a pair
$(H,\|\cdot\|)$, where $H\in\bbr^{m\times M}$ and \mbox{$\|\cdot\|$} is a
norm on $\bbr^M$, satisfies the condition $\bQ_{s,q}(\kappa)$
associated with the matrix $A$ and the r.s. $\cS$, if
%
\begin{equation}
\label{apeq1} \forall x\in\bbr^n\qquad  L_{s,q}(Bx)\leq
s^{{1}/{q}}\bigl\|H^TAx\bigr\|+\kappa s^{{1}/{q}-1}L_1(Bx).
\end{equation}

The following observation is evident:
%
\begin{observation}\label{apobserv1}
Given $A$ and an r.s. $\cS$, let $(H,\|\cdot\|)$ satisfy $\bQ
_{s,q}(\kappa)$. Then $(H, \|\cdot\|)$ satisfies $\bQ_{s,q'}(\kappa')$
for all $q'\in(1,q)$ and $\kappa'\geq\kappa$. Besides this, if $s'\leq
s$ is a positive integer, $((s/s')^{{1/q}}H,\|\cdot\|)$ satisfies
$\bQ_{s',q}((s'/s)^{1-{1}/{q}}\kappa)$. Furthermore, if
$(H,\|\cdot\|)$ satisfies $\bQ_{s,q}(\kappa)$, and $q'\ge q$, a
positive integer $s'\leq s$, and $\kappa'$ are such that $\kappa'
(s')^{{1}/{q'}-1}\ge\kappa s^{{1}/{q}-1}$, then
$({s^{{1}/{q}}(s')^{-{1}/{q'}}}H, \|\cdot\|)$ satisfies
$\bQ_{s',q'}(\kappa')$. In particular, when $s'\le s^{1-{1}/{ q}}$,
the fact that $(H,\|\cdot\|)$ satisfies $\bQ_{s,q}(\kappa)$ implies
that $(s^{1/q} H, \|\cdot\|)$ satisfies
$\bQ_{s',\infty}(\kappa)$.
\end{observation}

\subsection*{Relation to known conditions for the validity of sparse
$\ell_1$ recovery}
Note that whenever
\[
\cS=\bigl(B,n_1,\ldots,n_K,\|\cdot\|_{(1)},\ldots,\|\cdot
\|_{(K)}\bigr)
\]
is the standard r.s., the condition $\bQ_{s,q}(\kappa)$ reduces to
the condition $\bH_{s,q}(\kappa)$ introduced in \cite{JNnoisy}.
On the other hand, condition $\bQ_{s,p}(\kappa)$ is closely related
to other known conditions, introduced to study the properties of
recovery routines in the context of block-sparsity. Specifically,
consider an r.s. with $B=I_n$, and let us make the following observation:

Let $(H,\|\cdot\|_\infty)$ satisfy $\bQ_{s,q}(\kappa)$ and let
$\widehat{\lambda}$ be the maximum of the Euclidean norms of columns
in $H$. Then
%
\begin{equation}\label{sscond}
\forall x\in\bbr^n\qquad L_{s,q}(x)\leq\widehat{
\lambda}s^{{1/q}}\| Ax\|_2+\kappa s^{{1/ q}-1}L_1(x).
\end{equation}
Let us fix the r.s. $\cS_2=(I_n, n_1,\ldots, n_K, \|\cdot\|_2,\ldots,\|
\cdot\|_2)$.
Condition (\ref{sscond}) with $\kappa<1/2$ plays a crucial role in the
performance analysis\vadjust{\goodbreak} of the group-Lasso and Dantzig Selector. For
example, the error bounds for Lasso recovery obtained in \cite
{Lounicietal10} rely upon the Restricted Eigenvalue assumption $\RE
(s,\varkappa)$ as follows:
there exists $\varkappa>0$ such that
\[
L_2\bigl(x^s\bigr)\leq\frac{1}{\varkappa}\|Ax
\|_2\qquad\mbox{whenever } 3 L_1\bigl(x^s\bigr)
\geq L_1\bigl(x-x^s\bigr).
\]
In this case $L_{s,1}(x)\leq\sqrt{s}L_{s,2}(x)\leq\frac{\sqrt{s}}{
\varkappa}\|Ax\|_2$ whenever $4L_{s,1}(x)\geq L_1(x)$, so that
%
\begin{equation}
\label{Tsybakov} \forall x\in\bbr^n\qquad  L_{s,1}(x)\leq
\frac{s^{1/2}}{\varkappa}\|Ax\|_2+\frac{1}{4}L_1(x),
\end{equation}
which is exactly (\ref{sscond}) with $q=1$, $\kappa=1/4$ and
$\widehat
{\lambda}=(\varkappa\sqrt{s})^{-1}$ (observe that (\ref{Tsybakov}) is
nothing but the ``block version'' of the Compatibility condition from~\cite{BuhlmannvandeGeer09}).

Recall that a sensing matrix $A\in\bbr^{m\times n}$ satisfies the
Block Restricted Isometry Property $\BRIP(\delta,k)$ (see, e.g.,
\cite{EldarMishali09}) with $\delta\geq0$ and
a positive integer $k$ if for every $x\in\bbr^n$ with at most $k$
nonvanishing blocks one has
%
\begin{equation}
\label{apeq100bl} (1-\delta)\|x\|_2^2\leq
x^TA^TAx\leq(1+\delta)\|x\|_2^2.
\end{equation}

\begin{proposition}\label{appropRIPBlock} Let $A\in\bbr^{m\times n}$
satisfy $\BRIP(\delta,2s)$ for some $\delta<1$ and positive integer
$s$. Then:

\begin{longlist}
\item
The pair $ (H=\frac{s^{-1/2}}{\sqrt{1-\delta}}I_m,\|
\cdot\|_2 )$ satisfies the condition $\bQ_{s,2} (\frac
{\delta
}{1-\delta} )$ associated with $A$
and the r.s. $\cS_2$.

\item The pair $ (H=\frac{1}{1-\delta}A,L_\infty(\cdot
) )$ satisfies the condition $\bQ_{s,2} (\frac{\delta}{
1-\delta} )$ associated with $A$ and the r.s. $\cS_2$.
\end{longlist}
\end{proposition}
Our last observation here is as follows: let $(H,\|\cdot\|)$ satisfy
$\bQ_{s,q}(\kappa)$ for the r.s. given by $(B, n_1,\ldots, n_K, \|\cdot
\|_2,\ldots,\|\cdot\|_2)$, and let $d=\max_kn_k$. Then $(H,\|\cdot\|)$
satisfies $\bQ_{s,q}(\sqrt{d}\kappa)$ for the r.s. given by \mbox{$(B,
n_1,\ldots, n_K, \|\cdot\|_\infty,\ldots,\|\cdot\|_\infty)$}.

\section{\texorpdfstring{Accuracy bounds for $\ell_1$ block recovery routines}{Accuracy bounds for l1 block recovery routines}}
\label{secrecovery}

Throughout this section we fix an r.s. $\cS=(B,n_1,\ldots,n_K,\|\cdot\|
_{(1)},\ldots,\|\cdot\|_{(K)})$ and a sensing matrix~$A$.

\subsection{\texorpdfstring{Regular $\ell_1$ recovery}{Regular l1 recovery}}
We define the \textit{regular $\ell_1$ recovery} as
%
\begin{equation}
\label{apeq4} \widehat{x}_\reg(y)\in\Argmin_u \bigl
\{L_1(Bu)\dvtx \bigl\|H^T(Au-y)\bigr\|\leq\rho\bigr\},
\end{equation}
where the \textit{contrast matrix} $H\in\bbr^{m\times M}$, the norm $\|
\cdot\|$ and $\rho>0$ are parameters of the construction.
%
\begin{theorem}\label{aptheorem1} Let $s$ be a positive integer, $q\in
[1,\infty]$, $\kappa\in(0,1/2)$. Assume that the pair $(H,\|\cdot\|
)$ satisfies the condition $\bQ_{s,q}(\kappa)$ associated with $A$
and r.s. $\cS$, and let
%
\begin{equation}
\label{Xiofrho} \Xi=\Xi_{\rho,\cU}=\bigl\{\xi\dvtx  \bigl\|H^T(u+\xi)\bigr\|
\leq\rho\ \forall u\in\cU\bigr\}.\vadjust{\goodbreak}
\end{equation}
Then for all $x\in\bbr^n$, $u\in{\cal U}$ and $\xi\in\Xi$ one has
%
\begin{eqnarray}
\label{apeq6}
&&L_p\bigl(B\bigl[\widehat{x}_\reg(Ax+u+
\xi)-x\bigr]\bigr)\nonumber\\[-9pt]\\[-9pt]
&&\qquad\leq\frac{4(2s)^{{1}/{ p}}}{
1-2\kappa} \biggl[\rho+\frac{1}{2s}L_1
\bigl(Bx-[Bx]^s\bigr) \biggr],\qquad 1\leq p\leq q.\nonumber
\end{eqnarray}
\end{theorem}
%
The above result can be slightly strengthened by replacing the
assumption that $(H,\|\cdot\|)$ satisfies $\bQ_{s,q}(\kappa)$,
$\kappa<1/2$, with a weaker, by Observation~\ref{apobserv1},
assumption that $(H,\|\cdot\|)$ satisfies $\bQ_{s,1}(\varkappa)$
with $\varkappa<1/2$ and satisfies $\bQ_{s,q}(\kappa)$ with some
(perhaps large) $\kappa$:
%
\begin{theorem}\label{aptheorem101} Given $A$, r.s. $\cS$, integer
$s>0$, $q\in[1,\infty]$ and $\varepsilon\in(0,1)$, assume that $(H,\|
\cdot\|)$ satisfies the condition $\bQ_{s,1}(\varkappa)$ with
$\varkappa<1/2$ and the condition $\bQ_{s,q}(\kappa)$ with some
$\kappa\geq\varkappa$, and let $\Xi$ be given by (\ref{Xiofrho}).
Then for all $x\in\bbr^n$, $u\in{\cal U}$, $\xi\in\Xi$ and $p,
1\leq p\leq q$, it holds
%
\begin{eqnarray}
\label{apeq601}
&&
L_p\bigl(B\bigl[\widehat{x}_\reg(Ax+u+
\xi)-x\bigr]\bigr) \nonumber\\[-9pt]\\[-9pt]
&&\qquad\leq\frac{4(2s)^{{1}/{ p}}[1+\kappa-\varkappa
]^{{q(p-1)}/({
p(q-1)})}}{1-2\varkappa} \biggl[\rho+\frac{L_1(Bx-[Bx]^s)}{2s}
\biggr].\nonumber
\end{eqnarray}
\end{theorem}
%
\subsection{\texorpdfstring{Penalized $\ell_1$ recovery}{Penalized l1 recovery}}
The penalized $\ell_1$ recovery is
%
\begin{equation}
\label{apeq4pen} \widehat{x}_\pen(y)\in\Argmin_u \bigl
\{L_1(Bu)+\lambda\bigl\|H^T(Ax-y)\bigr\| \bigr\},
\end{equation}
where $H\in\bbr^{m\times M}$, $\|\cdot\|$ and a positive real
$\lambda$ are parameters of the construction.
%
\begin{theorem}\label{aptheorem2} Given $A$, r.s. $\cS$, integer $s$,
$q\in[1,\infty]$ and $\varepsilon\in(0,1)$, assume that $(H,\|\cdot\|
)$ satisfies the conditions $\bQ_{s,q}(\kappa)$ and $\bQ_{s,1}(\varkappa
)$ with $\varkappa<1/2$ and $\kappa\geq\varkappa$.

\begin{longlist}
\item
Let $\lambda\ge2s$.
Then for all $x\in\bbr^n$, $y\in\bbr^m$ it holds for $1\le p\le q$
%
\begin{eqnarray}
\label{apeq2601}
&&
L_p\bigl(B\bigl[\widehat{x}_\pen(y)-x
\bigr]\bigr)
\nonumber\\[-2pt]
&&\qquad\leq\frac{4\lambda^{{1}/{ p}}}{1-2\varkappa} \biggl[1+\frac
{\kappa
\lambda}{2s}-\varkappa
\biggr]^{{q(p-1)}/({ p(q-1)})} \\[-2pt]
&&\qquad\quad{}\times\biggl[\bigl\| H^T(Ax-y)\bigr\|+\frac{1}{2s}L_1
\bigl(Bx-[Bx]^s\bigr) \biggr].
\nonumber
\end{eqnarray}
In particular, with $\lambda=2s$ we have for $1\leq p\leq q$
%
\begin{eqnarray}
\label{apeq2601a}
&&
L_p\bigl(B\bigl[\widehat{x}_\pen(y)-x
\bigr]\bigr)
\nonumber\\
&&\qquad\leq\frac{4(2s)^{{1}/{
p}}}{1-2\varkappa} [1+\kappa-\varkappa]^{{q(p-1)}/({
p(q-1)})} \\
&&\qquad\quad{}\times\biggl[
\bigl\|H^T(Ax-y)\bigr\|+\frac{1}{
2s}L_1\bigl(Bx-[Bx]^s
\bigr) \biggr].
\nonumber
\end{eqnarray}

\item Let $\rho\geq0$ and $\Xi$ be given by (\ref
{Xiofrho}). Then for all $x\in\bbr^n$, $u\in{\cal U}$ and all $\xi
\in\Xi$ one has for $1\leq p\leq q$
%
\begin{eqnarray}
\label{apeq2626standard} \lambda&\geq&2s\quad\Rightarrow\quad L_p\bigl(B\bigl[
\widehat{x}_\pen(Ax+u+\xi)-x\bigr]\bigr)
\nonumber
\\
&&\qquad\hspace*{48pt}\leq\frac{4\lambda^{{1}/{ p}}}{1-2\varkappa} \biggl[1+\frac
{\kappa
\lambda}{2s}-\varkappa
\biggr]^{{q(p-1)}/({ p(q-1)})} \nonumber\\
&&\qquad\quad\hspace*{48pt}{}\times\biggl[\rho+\frac{1}{2s}L_1
\bigl(Bx-[Bx]^s\bigr) \biggr],
\\
\lambda&=&2s\quad\Rightarrow\quad L_p\bigl(B\bigl[\widehat{x}_\pen(Ax+u+
\xi)-x\bigr]\bigr)
\nonumber
\\
&&\qquad\hspace*{48pt}\leq\frac
{4(2s)^{{1}/{
p}}}{1-2\varkappa} [1+\kappa-\varkappa]^{{q(p-1)}/({ p(q-1)})} \nonumber\\
&&\qquad\quad\hspace*{48pt}{}\times\biggl
[\rho+
\frac
{1}{2s}L_1\bigl(Bx-[Bx]^s\bigr) \biggr].
\nonumber
\end{eqnarray}
\end{longlist}
\end{theorem}

\subsubsection*{Discussion} Let us compare the error bounds of the
regular and the penalized $\ell_1$ recoveries associated with the same
pair $(H,\|\cdot\|)$ satisfying the condition
$\bQ_{s,q}(\kappa)$ with $\kappa<1/2$. Given $\varepsilon\in(0,1)$, let
%
\begin{equation}
\label{rhoofH}\quad\rho_\varepsilon\bigl[H,\|\cdot\|\bigr]=\min\bigl\{\rho\dvtx  \Prob\bigl
\{
\xi\dvtx \bigl\|H^T(u+\xi)\bigr\|\leq\rho\ \forall u\in\cU\bigr\}\geq1-\varepsilon
\bigr\};
\end{equation}
this is nothing but the smallest $\rho$ such that
%
\begin{equation}
\label{probound} \Prob(\xi\in\Xi_{\rho,\varepsilon})\geq1-\varepsilon
\end{equation}
[see (\ref{Xiofrho})] and, thus, the smallest $\rho$ for which
the error bound (\ref{apeq6}) for the regular $\ell_1$ recovery holds
true with probability $1-\varepsilon$ (or at least the smallest $\rho$
for which the latter claim is supported by Theorem~\ref{aptheorem1}).
With $\rho=\rho_\varepsilon[H,\|\cdot\|]$, the regular $\ell_1$ recovery guarantees (and that is the best guarantee one can extract
from Theorem~\ref{aptheorem1}) that

\begin{quote}
(!) For some set $\Xi$, $\Prob\{\xi\in\Xi\}\geq1-\varepsilon$,
of ``good'' realizations of the random component $\xi$ of the
observation error, one has
%
\begin{eqnarray}
\label{theboundreads}
&&
L_p\bigl(B\bigl[\widehat{x}(Ax+u+\xi)-x\bigr]
\bigr)\nonumber\\[-8pt]\\[-8pt]
&&\qquad\leq\frac{4(2s)^{{1}/{ p}}}{
1-2\kappa} \biggl[\rho_\varepsilon\bigl[H,\|\cdot\|\bigr]+
\frac{L_1(Bx-[Bx]^s)}{
2s} \biggr],\qquad 1\leq p\leq q,\nonumber
\end{eqnarray}
whenever $x\in\bbr^n,u\in\cU$, and $\xi\in\Xi$.
\end{quote}

\noindent The error bound (\ref{apeq2601a}) [where we can safely set $\varkappa
=\kappa$, since $\bQ_{s,q}(\kappa)$ implies $\bQ_{s,1}(\kappa)$]
says that \textit{(!) holds true for the penalized $\ell_1$ recovery
with $\lambda=2s$.}
The latter observation suggests that the penalized $\ell_1$ recovery
associated with $(H,\|\cdot\|)$ and $\lambda=2s$ is better than its
regular counterpart, the reason being twofold.\vadjust{\goodbreak} First, in order to
ensure (!) with the regular recovery, the ``built in'' parameter $\rho
$ of this recovery should be set to $\rho_\varepsilon[H,\|\cdot\|]$,
and the latter quantity is not always easy to identify. In contrast to
this, the construction of the penalized $\ell_1$ recovery is
completely independent of a priori assumptions on the structure of
observation errors, while automatically ensuring (!) for the error
model we use. Second, and more importantly, for the penalized recovery
the bound (\ref{theboundreads}) is no more than the ``worst, with
confidence $1-
\varepsilon$, case,'' and the typical values of the quantity $\|H^T(u+\xi
)\|$ which indeed participates in the error bound (\ref{apeq2601}) are
essentially smaller than $\rho_\varepsilon[H,\|\cdot\|]$. Our numerical
experience fully supports the above suggestion: the difference in
observed performance of the two routines in question, although not
dramatic, is definitely in favor of the penalized recovery. The only
potential disadvantage of the latter routine is that the penalty
parameter $\lambda$ should be tuned to the level $s$ of sparsity we
aim at, while the regular recovery is free of any guess of this type.
Of course, the ``tuning'' is rather loose---all we need (and
experiments show that we indeed need this) is the relation $\lambda
\geq2s$, so that a rough upper bound on $s$ will do; note, however,
that the bound (\ref{apeq2601}) deteriorates as $\lambda$ grows.

\section{\texorpdfstring{Tractability of condition $\bQ_{s,\infty}(\kappa)$, $\ell_\infty$-norm of the blocks}
{Tractability of condition Q s,infinity(kappa), l infinity-norm of the blocks}}\label{secconditionproperties}

We have seen in Section~\ref{secrecovery} that given a sensing matrix
$A$ and an r.s. $\cS=(B,n_1,\ldots,n_K,\break\|\cdot\|_{(1)},\ldots,\|\cdot\|
_{(K)})$ such that the associated conditions $\bQ_{s,q}(\kappa)$ are
satisfiable, we can validate the $\ell_1$ recovery of nearly
$s$-block-sparse signals, specifically, we can point out $\ell_1$-type
recoveries with controlled (and small, provided so are the observation
error and the deviation of the signal from an $s$-block-sparse one).
The bad news here is that, in general, condition $\bQ_{s,q}(\kappa)$,
as well as other conditions for the validity of $\ell_1$ recovery,
like Block RE/RIP, cannot be verified efficiently. The latter means
that given a sensing matrix $A$ and a r.s. $\cS$, it is difficult to
verify that a given candidate pair $(H,\|\cdot\|)$ satisfies condition
$\bQ_{s,q}(\kappa)$ associated with $A$ and $\cS$.
Fortunately, one can construct ``tractable approximations'' of
condition $\bQ_{s,q}(\kappa)$, that is, verifiable sufficient
conditions for the validity of $\bQ_{s,q}(\kappa)$. The first good
news is that when all $\|\cdot\|_{(k)}$ are the uniform norms $\|\cdot
\|_\infty$ and, in addition, $q=\infty$ [which, by Observation \ref
{apobserv1}, corresponds to the strongest among the conditions $\bQ
_{s,q}(\kappa)$ and ensures the validity of (\ref{apeq6}) and (\ref
{apeq2601}) in the largest possible range $1\leq p\leq\infty$ of
values of $p$], the condition $\bQ_{s,q}(\kappa)$ becomes ``fully
computationally tractable.'' We intend to demonstrate also that the
condition $\bQ_{s,\infty}(\kappa)$ is in fact necessary for the risk
bounds of the form (\ref{apeq6})--(\ref{apeq2626standard}) to be
valid when $p=\infty$.

\subsection{\texorpdfstring{Condition $\bQ_{s,\infty}(\kappa)$: Tractability and the optimal choice of the contrast matrix $H$}
{Condition Q s,infinity(kappa): Tractability and the optimal choice of the contrast matrix H}}
\label{secconditiontractability}

\subsubsection*{Notation} In the sequel, given $r,\theta\in[1,\infty
]$ and a matrix $M$, we denote by $\|M\|_{r,\theta}$ the norm of the
linear operator $u\mapsto Mu$ induced\vadjust{\goodbreak} by the norms $\|\cdot\|_r$
and
$\|\cdot\|_\theta$ on the argument and the image spaces:
\[
\|M\|_{r,\theta}=\max_{u\dvtx \|u\|_r\leq1}\|Mu\|_\theta.
\]
We denote by $\|M\|_{(\ell,k)}$ the norm of the linear mapping
$u\mapsto Mu\dvtx \bbr^{n_\ell}\to\bbr^{n_k}$ induced by the norms $\|
\cdot\|_{(\ell)}$, $\|\cdot\|_{(k)}$ on the argument and on the
image spaces.
Further, $\Row_k[M]$ stands for the transpose of the $k$th row of $M$
and $\Col_k[M]$ stands for $k$th column of $M$. Finally, $\|u\|_{s,q}$
is the $\ell_q$-norm of the vector obtained from a vector $u\in\bbr^k$
by zeroing all but the $s$ largest in magnitude entries in $u$.

\subsubsection*{Main result}
Consider r.s. $\cS_\infty=(B,n_1,\ldots,n_K, \|\cdot\|_{\infty},\ldots,\|
\cdot\|_{\infty})$. We claim that in this case the condition $\bQ
_{s,\infty}(\kappa)$ becomes fully tractable. Specifically, we have
the following.
%
\begin{proposition}\label{appropTract} Let a matrix $A\in\bbr^{m\times
n}$, the r.s. $\cS_\infty$, a positive integer $s$ and
reals $\kappa>0$, $\varepsilon\in(0,1)$ be given.

\begin{longlist}
\item
Assume that a triple $(H,\|\cdot\|,\rho)$, where $H\in\bbr
^{m\times M}$, $\|\cdot\|$ is a norm on $\bbr^M$, and $\rho\geq0$,
is such that

\begin{quote}
\textup{(!)}
$(H,\|\cdot\|)$ satisfies $\bQ_{s,\infty}(\kappa)$, and the set
$\Xi=\{\xi\dvtx \|H^T[u+\xi]\|\leq\rho\ \forall u\in\cU\}$
satisfies $\Prob(\xi\in\Xi)\geq1-\varepsilon$.
\end{quote}

\noindent Given $H$, $\|\cdot\|$, $\rho$, one can find efficiently
$N=n_1+\cdots+n_K$ vectors $h^1,\ldots,h^N$ in $\bbr^m$ and $N\times N$
block matrix $V=[V^{k\ell}]_{k,\ell=1}^K$ (the blocks $V^{k\ell}$ of
$V$ are $n_k\times n_\ell$ matrices) such that
%
\begin{eqnarray}
\label{neweq333}\quad &&\mbox{\textup{(a)}}\quad B=VB+\bigl[h^1,\ldots,h^N
\bigr]^TA,
\nonumber
\\
&&\mbox{\textup{(b)}}\quad\bigl\|V^{k\ell}\bigr\|_{\infty,\infty}\leq s^{-1}\kappa\qquad
\forall
k,\ell\leq K,
\\
&&\mbox{\textup{(c)}}\quad\Prob_\xi\Bigl(\Xi^+:=\Bigl\{\xi\dvtx  \max_{u\in\cU
}u^Th^i+\bigl|
\xi^Th^i\bigr|\leq\rho, 1\leq i\leq N\Bigr\} \Bigr)\geq1-
\varepsilon
\nonumber
\end{eqnarray}
(note that the matrix norm $\|A\|_{\infty, \infty}=\max_j\|\Row_j[A]\|
_1$ is simply the maximum $\ell_1$-norm of the rows of $A$).

\item Whenever vectors $h^1,\ldots,h^N\in\bbr^m$ and a matrix
$V=[V^{k\ell}]_{k,\ell=1}^K$ with $n_k\times n_\ell$ blocks
$V^{k\ell}$ satisfy (\ref{neweq333}), the $m\times N$ matrix
$\widehat{H}=[h^1,\ldots,h^N]$, the norm $\|\cdot\|_\infty$ on $\bbr^N$
and $\rho$ form a triple satisfying $(!)$.
\end{longlist}
\end{proposition}

\subsubsection*{Discussion} Let a sensing matrix $A\in\bbr^{m\times
n}$ and a r.s. $\cS_\infty$ be given, along with a positive integer
$s$, an uncertainty set $\cU$, a distribution $P$ of $\xi$ and
$\varepsilon\in(0,1)$. Theorems~\ref{aptheorem1} and~\ref{aptheorem2} say
that if a triple $(H,\|\cdot\|,\rho)$ is such that $(H,\|\cdot\|)$
satisfies $\bQ _{s,\infty}(\kappa)$ with $\kappa<1/2$ and $H,\rho$ are
such that the set $\Xi$ given by (\ref{Xiofrho}) satisfies
(\ref{probound}), then for the regular $\ell_1$ recovery associated
with $(H,\|\cdot\|,\rho)$ and for the penalized $\ell_1$ recovery
associated with $(H,\| \cdot\|)$ and $\lambda=2s$, the following holds:
%
\begin{eqnarray}
\label{apeq60}
&&\forall\bigl(x\in\bbr^n,u\in\cU,\xi\in\Xi\bigr)\qquad
\nonumber\\
&&\qquad L_p\bigl(B\bigl[\widehat{x}(Ax+u+\xi)-x\bigr]\bigr)\leq
\frac{4(2s)^{{1}/{ p}}}{
1-2\kappa} \biggl[\rho+\frac{1}{2s}L_1
\bigl(Bx-[Bx]^s\bigr) \biggr],\\
&&\eqntext{1\leq p\leq\infty.}
\end{eqnarray}
Proposition~\ref{appropTract} states that when applying this result,
we lose nothing by restricting ourselves with triples
$H=[h^1,\ldots,h^N]\in\bbr^{m\times N}$, $N=n_1+\cdots+n_K$, $\|\cdot\|
=L_\infty(\cdot)$, $\rho\geq0$ which can be augmented by an
appropriately chosen $N\times N$ matrix $V$ to satisfy relations (\ref
{neweq333}). In the rest of this discussion, it is assumed that we are
speaking about triples $(H,\|\cdot\|,\rho)$ satisfying the just
defined restrictions.

The bound (\ref{apeq60}) is completely determined by two
parameters---$\kappa$ (which should be $<1/2$) and $\rho$; the smaller
are these parameters, the better are the bounds. In what follows we
address the issue of efficient synthesis of matrices $H$ with ``as good
as possible'' values of $\kappa$ and $\rho$.

Observe first that $H=[h^1,\ldots,h^N]$ and $\kappa$ should admit an
extension by a matrix $V$ to a solution of the system of convex
constraints (\ref{neweq333})(a), (\ref{neweq333})(b). In the case of
$\xi\equiv0$ the best choice of $\rho$, given $H$, is
\[
\rho=\max_i\mu_{\cU}\bigl(h^i\bigr) \qquad\mbox{where } \mu_{\cU
}(h)=\max_{u\in\cU} u^Th.
\]
Consequently, in this case the ``achievable pairs'' $\rho$, $\kappa$
form a computationally tractable convex set
\begin{eqnarray*}
G_s&=& \biggl\{(\kappa,\rho)\dvtx  \exists H=\bigl[h^1,\ldots,h^N
\bigr]\in\bbr^{m\times
N},\\
&&\hspace*{7pt}V=\bigl[V^{k\ell}\in\bbr^{n_k\times n_\ell}
\bigr]_{k,\ell=1}^K\dvtx
B=VB+H^TA, \\
&&\hspace*{27.5pt}\bigl\|V^{k\ell}\bigr\|_{\infty,\infty}\leq \frac{\kappa}{
s},\mu_{\cU}\bigl(h^i\bigr)\leq\rho,1\leq i\leq N \biggr\}.
\end{eqnarray*}
When $\xi$ does not vanish, the situation is complicated by the
necessity to maintain the validity of the restriction
%
\begin{eqnarray}
\label{apeq70}
\Prob\bigl(\xi\in\Xi^+\bigr):\!&=&\Prob\bigl\{\xi\dvtx
\mu_{\cU}\bigl(h^i\bigr)+\bigl|\xi^Th^i\bigr|
\leq\rho, 1\leq i\leq N \bigr\}\nonumber\\[-8pt]\\[-8pt]
&\geq&1-\varepsilon,\nonumber
\end{eqnarray}
which is a chance constraint in variables $h^1,\ldots,h^N,\rho$ and as
such can be ``computationally intractable.'' Let us consider the
``standard'' case of Gaussian zero mean noise $\xi$, that is, assume
that $\xi=D\eta$ with $\eta\sim\cN(0,I_m)$ and known $D\in\bbr^{m\times
m}$. Then
(\ref{apeq70}) implies that
\[
\rho\geq\max_i \biggl[\mu_{\cU}\bigl(h^i
\bigr)+\ErfInv\biggl(\frac{\varepsilon}{2} \biggr)\bigl\|D^Th^i
\bigr\|_2 \biggr].
\]
On the other hand, (\ref{apeq70}) is clearly implied by
\[
\rho\geq\max_i \biggl[ \mu_{\cU}\bigl(h^i
\bigr)+\ErfInv\biggl(\frac{\varepsilon}{2N} \biggr)\bigl\|D^Th^i
\bigr\|_2 \biggr].
\]
Ignoring the ``gap'' between $\ErfInv(\frac{\varepsilon}{2}
)$ and $\ErfInv(\frac{\varepsilon}{2N} )$, we can safely
model the restriction (\ref{apeq70}) by the system of convex constraints
%
\begin{equation}
\label{apeq71} \mu_{\cU}\bigl(h^i\bigr)+\ErfInv\biggl(
\frac{\varepsilon}{2N} \biggr)\bigl\|D^Th^i\bigr\|_2\leq
\rho,\qquad 1\leq i\leq N.
\end{equation}
Thus, the set $G_s$ of admissible $\kappa,\rho$ can be safely
approximated by the computationally tractable convex set
%
\begin{eqnarray}
\label{apeq72} && G_s^*=\biggl\{(\kappa,\rho)\dvtx \exists\bigl[
H=\bigl[h^1,\ldots,h^N\bigr]\in
\bbr^{m\times N},\nonumber\\
&&\hspace*{53.5pt}\hspace*{33.5pt} V=\bigl[V^{k\ell}\in\bbr^{n_k\times n_\ell}
\bigr]_{k,\ell=1}^K \bigr]\dvtx\nonumber\\[-8pt]\\[-8pt]
&&\hspace*{33.5pt}
\biggl\{ B=BV+H^TA, \bigl\|V^{k\ell}\bigr\|_{\infty,\infty}\leq
\frac{\kappa}{ s}, 1\leq k,\ell\leq K,
\nonumber
\\
&&\hspace*{41pt} \max_{u\in\cU}u^Th^i+\ErfInv\biggl(
\frac{\varepsilon}{
2N} \biggr)\bigl\|D^Th^i\bigr\|_2\leq
\rho, 1\leq i\leq N\biggr\}\biggr\}.\nonumber
\end{eqnarray}

\subsection{\texorpdfstring{Condition $\bQ_{s,\infty}(\kappa)$: Necessity}
{Condition Q s,infinity(kappa): Necessity}}
\label{secconditionnecessity}

In this section, as above, we assume that all norms $\|\cdot\|_{(k)}$
in the r.s. $\cS_\infty$ are $\ell_\infty$-norms; we assume, in
addition, that $\xi$ is a zero mean Gaussian noise: $\xi=D\eta$ with
$\eta\sim\cN(0,I_m)$ and known $D\in\bbr^{m\times m}$. From the
above discussion we know that if, for some $\kappa<1/2$ and $\rho>0$,
there exist $H=[h^1,\ldots,h^N]\in\bbr^{m\times N}$ and $V=[V^{k\ell
}\in\bbr^{n_k\times n_\ell}]_{k,\ell=1}^K$ satisfying (\ref
{neweq333}), then regular and penalized $\ell_1$ recoveries with
appropriate choice of parameters ensure that
%
\begin{eqnarray}
\label{apeq80}\quad
&&\forall\bigl(x\in\bbr^n,u\in\cU\bigr)
\nonumber\\
&&\qquad\Prob_\xi\biggl(\bigl\|B\bigl[x-\widehat{x}(Ax+u+\xi)\bigr]
\bigr\|_\infty\leq\frac{4}{
1-2\kappa} \biggl[\rho+\frac{L_1(Bx-[Bx]^s)}{2s} \biggr]
\biggr)\\
&&\qquad\qquad\geq1-\varepsilon.
\nonumber
\end{eqnarray}
We are about to demonstrate that this implication can be ``nearly
inverted'':

\begin{proposition}\label{appropNec}
Let a sensing matrix $A$, an r.s. $\cS_\infty$ with $\|\cdot\|_{(k)}=\|
\cdot\|_\infty$, $1\leq k\leq K$, an uncertainty set $\cU$,
and reals $\kappa>0$, $\varepsilon\in(0,1/2)$ be given. Suppose that
the observation error ``is present,'' specifically, that for every
$r>0$, the set $\{u+ De\dvtx u\in\cU,\|e\|_2\leq r\}$ contains a
neighborhood of the origin.\looseness=-1

Given a positive integer $S$, assume that there exists a recovering
routine $\widehat{x}$ satisfying an error bound of the form
(\ref{apeq80}), specifically, such that for all $x\in\bbr^n,u\in\cU$,
%
\begin{equation}
\label{apeq81} \Prob_\xi\bigl( \bigl\|B\bigl[x-\widehat{x}(Ax+u+\xi)
\bigr]\bigr\|_\infty\leq\alpha+S^{-1}L_1
\bigl(Bx-[Bx]^S\bigr) \bigr)\geq1-\varepsilon\hspace*{-32pt}\vadjust{\goodbreak}
\end{equation}
for some $\alpha>0$. Then there exist $H=[h^1,\ldots,h^N]\in\bbr^{m\times
N}$ and $V=[V^{k\ell}\in\bbr^{n_k\times n_\ell}]_{k,\ell
=1}^K$ satisfying
%
\begin{eqnarray}
\label{apeq88} &&\mbox{\textup{(a)}}\quad B=VB+H^TA,
\nonumber
\\
&&\mbox{\textup{(b)}}\quad\bigl\|V^{k\ell}\bigr\|_{\infty,\infty}\leq{2 S^{-1}}\qquad \forall k,\ell
\leq K,
\\
&&\mbox{\textup{(c)}}\quad\mbox{with }\rho:=\max_{1\leq i\leq N} \biggl[\max_{u\in\cU}u^Th^i+{
\ErfInv} \biggl(\frac{\varepsilon}{
2N} \biggr)\bigl\|D^Th^i
\bigr\|_2 \biggr],
\nonumber
\end{eqnarray}
one has $\rho\leq2\alpha$ when $D=0$, $\rho\leq2\alpha\frac
{\sErfInv
({\varepsilon}/({2N}))}{\sErfInv(\varepsilon)}$ when $D\neq0$,
and for $\xi=D\eta, \eta\sim\cN(0,I_m)$ one has
\[
\Prob_\xi\Bigl(\Xi^+:=\Bigl\{\xi\dvtx  \max_{u\in\cU}u^Th^i+\bigl|
\xi^Th^i\bigr|\leq\rho, 1\leq i\leq N\Bigr\} \Bigr)\geq1-
\varepsilon.
\]
\end{proposition}
In other words (see Proposition~\ref{appropTract}), $(H,L_\infty
(\cdot))$ satisfies $\bQ_{s,\infty}(\kappa)$ for $s$ ``nearly as
large as $S$,'' namely, $s\leq\frac{\kappa}{2} S$, and
$H=[h^1,\ldots,h^k],\rho$
satisfy conditions (\ref{apeq71})
with $\rho$ being ``nearly $\alpha$,'' namely, $\rho\leq2\alpha$
in the case of $D=0$ and $\rho\leq2\frac{\sErfInv({\varepsilon}/({
2N}))}{\sErfInv(\varepsilon)}$ when $D\neq0$.
In particular, under the premise of Proposition~\ref{appropNec}, the
contrast optimization procedure of Section \ref
{secconditiontractability} supplies
the matrix $H$ such that the corresponding regular or penalized
recovery $ \widehat{x}(\cdot)$ for all $s\le\frac{S}{8}$ satisfies
\[
\Prob_\xi\biggl\{\bigl\|B\bigl[x-\widehat{x}(y)\bigr]\bigr\|_\infty
\leq4 \biggl[4\frac{\sErfInv({\varepsilon}/({2N}))}{\sErfInv(\varepsilon
)}\alpha+s^{-1}L_1
\bigl(Bx-[Bx]^s\bigr) \biggr] \biggr\}\geq1-\varepsilon.
\]
%
\section{\texorpdfstring{Tractable approximations of $\bQ_{s,q}(\kappa)$}
{Tractable approximations of Q s,q(kappa)}}
\label{secconditionsufficiency}

Aside from the important case of $q=\infty$, $\|\cdot\|_{(k)}=\|\cdot
\|_\infty$ considered in Sections~\ref{secconditiontractability}
and~\ref{secconditionnecessity}, condition $\bQ_{s,q}(\kappa)$
``as it is'' seems to be computationally intractable: unless $s=O(1)$,
it is unknown how to check efficiently that a given pair $(H,\|\cdot\|
)$ satisfies this condition, not speaking about synthesis of a pair
satisfying this condition and resulting in the best possible error
bound (\ref{apeq6}), (\ref{apeq2601}) for regular and penalized $\ell
_1$-recoveries.
We are about to present \textit{verifiable sufficient conditions} for the
validity of $\bQ_{s,q}(\kappa)$ which may become an interesting
alternative for condition $\bQ_{s,q}(\kappa)$ for that purposes.

\subsection{\texorpdfstring{Sufficient condition for $\bQ_{s,q}(\kappa)$}{Sufficient condition for Q s,q(kappa)}}
\label{secconditionsuff}

\begin{proposition}\label{propnec} Suppose that a sensing matrix $A$,
an r.s. $\cS=(B,\break n_1,\ldots, n_K,\|\cdot\|_{(1)},\ldots,\|\cdot\|_{(K)})$,
and $\kappa\geq0$ are given.

Let $N=n_1+\cdots+n_K$, and let $N\times N$ matrix $V=[V^{k\ell}]_{k,\ell
=1}^K$ ($V^{k\ell}$ are $n_k\times n_\ell$) and $m\times N$ matrix
$H$ satisfy the relation
%
\begin{equation}
\label{equationagain} B=VB+ H^TA.
\end{equation}
Let us denote
\[
\nu^*_{s,q}(V) = \max_{1\leq\ell\leq K} \max_{w^\ell\in\bbr
^{n_\ell}\dvtx \|w^\ell\|_{(\ell)}\leq1}L_{s,q}
\bigl(\bigl[V^{1\ell}w^\ell;\ldots;V^{K\ell}w^\ell
\bigr] \bigr).
\]
Then for all $s\leq K$ and all $q\in[1,\infty]$, we have
%
\begin{equation}
\label{onehasjanuary} L_{s,q}(Bx)\leq s^{{1/q}}L_\infty
\bigl(H^TAx\bigr) + \nu^*_{s,q}(V) L_1(Bx) \qquad \forall x\in\bR^n.
\end{equation}
\end{proposition}
The result of Proposition~\ref{propnec} is a step toward developing a
verifiable sufficient condition for the validity of $\bQ_{s,q}$. To
get such a condition, we need an efficiently computable upper bound of
the quantity $\nu^*_{s,q}$.
In particular, if for a given positive integer $s\leq K$ and a real
$q\in[1,\infty]$ there exist an upper bounding function $\nu_{s,q}(V)$
such that
%
\begin{equation}
\label{nustar}\nu_{s,q}(\cdot) \mbox{ is convex}\quad \mbox{and}\quad
\nu_{s,q}(V)\geq\nu^*_{s,q}(V)\qquad \forall V
\end{equation}
and a matrix $V$ such that
%
\begin{equation}
\label{goodcase} \nu_{s,q}(V)\leq s^{{1/ q}-1}\kappa,
\end{equation}
then the pair $(H,L_\infty(\cdot))$ satisfies $\bQ_{s,q}(\kappa)$.
An important example of the upper bound for $\nu^*_{s,q}(V)$ which
satisfies (\ref{goodcase}) is provided in the following statement.
%
\begin{proposition}\label{omegamatrix}
Let $\Omega$ be a $K\times K$ matrix with entries $[\Omega]_{k,\ell
}=\break\|V^{k\ell}\|_{(\ell,k)}$, $1\le k,\ell\le K$. Then
%
\begin{equation}
\label{defnuV} \widehat{\nu}_{s,q}(V):=\max_{1\leq k\leq K}\bigl\|
\Col_k[\Omega]\bigr\|_{s,q}\geq\nu^*_{s,q}(V)\qquad
\forall V
\end{equation}
[note that the inequality in (\ref{defnuV}) becomes equality when
either $q=\infty$ or $s=1$],
so that the condition
%
\begin{equation}
\label{sothatthe} \widehat{\nu}_{s,q}(V)\leq s^{{1/ q}-1}\kappa
\end{equation}
taken along with (\ref{equationagain}) is sufficient for
$(H,L_\infty(\cdot))$ to satisfy $\bQ_{s,q}(\kappa)$.
\end{proposition}
When all $\|\cdot\|_{(k)}$ are the $\ell_\infty$-norms and $q=\infty
$, the results of Propositions~\ref{propnec} and~\ref{omegamatrix}
recover Proposition~\ref{appropTract}. In the general case, they
suggest a way to synthesize matrices $H\in\bbr^{m\times N}$ which,
augmented by the norm $\|\cdot\|=L_{\infty}(\cdot)$, provably
satisfies the condition $\bQ_{s,q}(\kappa)$, along with a certificate
$V$ for this fact. Namely, $H$~and $V$ should satisfy the system of
linear equations (\ref{equationagain}) and, in addition, (\ref
{goodcase}) should hold for $V$ with $\nu_{s,q}(\cdot)$ satisfying
(\ref{nustar}). Further, for such a $\nu_{s,q}(\cdot)$, (\ref
{goodcase}) is a system of convex constraints on~$V$. Whenever these
constraints are efficiently computable, we get a computationally
tractable sufficient condition for $(H,L_\infty(\cdot))$ to satisfy
$\bQ_{s,q}(\kappa)$---a condition which is expressed by an explicit
system of efficiently computable convex constraints (\ref
{equationagain}), (\ref{goodcase}) on $H$ and additional matrix
variable $V$.

\subsection{Tractable sufficient conditions and contrast
optimization}\label{suffconddisc}
The quantity $\widehat{\nu}_{s,q}(\cdot)$ is the simplest choice of
$\nu_{s,q}(\cdot)$ satisfying (\ref{nustar}).
In this case, efficient computability of the constraints (\ref
{goodcase}) is the same as efficient computability of norms $\|\cdot\|
_{(k,\ell)}$. Assuming that $\|\cdot\|_{(k)}=\|\cdot\|_{r_k}$ for
every $k$ in the r.s. $\cS$, the computability issue becomes the one
of efficient computation of the norms $\|\cdot\|_{r_\ell,r_k}$. The
norm $\|\cdot\|_{r,\theta}$ is known to be generically efficiently
computable in only three cases:
\begin{longlist}[(3)]
\item[(1)]$\theta=\infty$, where $\|M\|_{r,\infty}=\|M^T\|_{1,{r}/({
r-1})}={\max_i}\|\Row_i^T(M)\|_{{r}/({ r-1})}$;
\item[(2)]$r=1$, where $\|M\|_{1,\theta}={\max_j}\|\Col_j[M]\|_\theta$;
\item[(3)]$r=\theta=2$, where $\|M\|_{2,2}=\sigma_{\max}(M)$ is the
spectral norm of $M$.
\end{longlist}
Assuming for the sake of simplicity that in our r.s. $\|\cdot\|_{(k)}$
are $r$-norms with common value of $r$, let us look at three
``tractable cases'' as specified by the above discussion---those of
$r=\infty$, $r=1$ and $r=2$. In these cases, candidate contrast
matrices $H$ are $m\times N$, the associated norm $\|\cdot\|$ is
$L_\infty(\cdot)$, and our sufficient condition for $H$ to be
\textit{good} [i.e., for $(H,L_\infty(\cdot))$ to satisfy
$\bQ_{s,q}(\kappa )$ with given $\kappa<1/2$ and $q$] becomes a system
$\bS=\bS_{\kappa,q}$ of explicit efficiently computable convex
constraints on $H$ and additional matrix variable $V\in\bbr^{N\times
N}$, implying that the set $\cH$ of good $H$ is convex and
computationally tractable, so that we can minimize efficiently over
$\cH$ any convex and efficiently computable function. In our context, a
natural way to use $\bS$ is to optimize over $H\in\cH$ the error bound
(\ref {theboundreads}) or, which is the same, to minimize over $\cH$
the function $\rho(H)=\rho_\varepsilon[H,L_\infty(\cdot)]$; see
(\ref{rhoofH}), where $\varepsilon<1$ is a given tolerance. Taken
literally, this problem still can be difficult, since the function
$\rho(H)$ is not necessarily convex and can be difficult to compute
even in the convex case. To overcome this difficulty, we again can use
a \textit{verifiable sufficient condition} for the relation $\rho
(H)\leq\rho$, that is, a system $\bT=\bT_\varepsilon$ of explicit
efficiently computable convex constraints on variables $H$ and $\rho$
(and, perhaps some slack variables $\zeta$) such that $\rho(H)\leq
\rho$ for the $(H,\rho)$-component of every feasible solution of $\bT
$. With this approach, the design of the best, as allowed by $\bS$ and
$\bT$, contrast matrix $H$ reduces to solving a convex optimization
problem with efficiently computable constraints in variables $H,V,\rho
$, specifically, the problem
%
\begin{equation}
\label{theconvexproblem} \min_{\rho,H,V,\zeta} \{\rho\dvtx  \mbox{$H,V$
satisfy $
\bS$; $H,\rho,\zeta$ satisfy $\bT$} \}.
\end{equation}
In the rest of this section we present explicitly the systems $\bS$
and $\bT$ for the three tractable cases we are interested in, assuming
the following model of observation errors:
\[
\cU=\bigl\{u=Ev\dvtx  \|v\|_2\leq1\bigr\};\qquad \xi=D\eta, \eta\sim\cN(0,I_m),
\]
where $E,D\in\bbr^{m\times m}$ are given.

We use the following notation: the $m\times N$ matrix $H$ is
partitioned into $m\times n_k$ blocks $H[k]$, $1\leq k\leq K$,
according to the block structure of the representation vectors; the
$t$th column in $H[k]$ is denoted $h^{kt}\in\bbr^m$, \mbox{$1\leq t\leq n_k$}.

For derivations of the results to follow, see Section
A.7 of the supplementary article \cite{JKNP-suppl}.

\subsubsection*{The case of $r=\infty$} The case of $q=\infty$ was
considered in full details in Section~\ref{secconditiontractability}.
When $q\leq\infty$, one has
%
\begin{eqnarray}
\label{caseinf} && \bS_{\kappa,q}\dvtx \cases{ %
B=VB+H^TA,
\vspace*{2pt}\cr
\displaystyle \Omega_{k\ell}:=\bigl\|V^{k\ell}\bigr\|_{\infty,\infty}=
\max_{1\leq t\leq
n_k}\bigl\|\Row_t\bigl[V^{k\ell}\bigr]
\bigr\|_1, &\quad $1\leq k,\ell\leq K$,
\vspace*{2pt}\cr
\bigl\|\Col_\ell[\Omega]
\bigr\|_{s,q}\leq s^{{1/ q}-1}\kappa, &\quad $1\leq\ell\leq K$,}
\hspace*{-32pt}\\
&& \bT_\varepsilon\dvtx \ErfInv\biggl(\frac{\varepsilon}{2N}\biggr)
\bigl\|D^Th^{kt}\bigr\|_2+\bigl\| E^Th^{kt}
\bigr\|_2\leq\rho,\nonumber\\
&&\eqntext{1\leq t\leq n_k, 1\leq k\leq K.}
\end{eqnarray}

\subsubsection*{The case of $r=2$} Here
%
\begin{eqnarray}
\label{casetwo} &&\bS_{\kappa,q}\dvtx
\cases{B=VB+H^TA,\vspace*{2pt}\cr
\Omega_{k\ell}:=\bigl\|V^{k\ell}\bigr\|_{2,2}=\sigma_{\max}
\bigl(V^{k\ell}\bigr), & \quad $1\leq k,\ell\leq K$,
\vspace*{2pt}\cr
\bigl\|\Col_\ell[\Omega]\bigr\|_{s,q}\leq s^{{1/ q}-1}\kappa, &\quad
$1 \leq\ell\leq K$,}
\nonumber\\
&&\bT_\varepsilon\dvtx \exists\bigl\{W_k\in\bbs^m,
\alpha_k,\beta_k,\gamma_k\in\bbr\bigr
\}_{k=1}^K\dvtx  \nonumber\\[-8pt]\\[-8pt]
&&\qquad\cases{\sigma_{\max}\bigl(E^TH[k]
\bigr)+\alpha_k\leq\rho,\vspace*{2pt}\cr
\lleft[\matrix{W_k&D^TH[k]
\vspace*{1pt}\cr
H^T[k]D&\alpha_kI_{n_k}} \rright]
\succeq0,\vspace*{3pt}\cr
\bigl\|\lambda(W_k)\bigr\|_\infty\leq\beta_k,\qquad \bigl\|
\lambda(W_k)\bigr\|_2\leq\gamma_k,\vspace*{2pt}\cr
\Tr(W_k)+2 \bigl[\delta\beta_k+\sqrt{
\delta^2\beta_k^2+2\delta
\gamma_k^2} \bigr]\leq\alpha_k,}\qquad
1\leq k\leq K,
\nonumber
\\
&&\delta:=\ln(K/\varepsilon),\nonumber
\end{eqnarray}
where $\bbs^m$ is the space of $m\times m$ symmetric matrices, and
$\lambda(W)$ is the vector of eigenvalues of $W\in\bbs^m$.

\subsubsection*{The case of $r=1$} Here
%
\begin{eqnarray}
\label{caseone}\qquad\quad &&\bS_{\kappa,q}\dvtx
\cases{B=VB+H^TA,\vspace*{2pt}\cr
\displaystyle \Omega_{k\ell}:=\bigl\|V^{k\ell}\bigr\|_{1,1}=
\max_{1\leq t\leq n_\ell}\bigl\| \Col_t\bigl[V^{k\ell}\bigr]
\bigr\|_1, &\quad $1\leq k,\ell\leq K$,\vspace*{2pt}\cr
\bigl\|\Col_\ell[\Omega]\bigr\|_{s,q}\leq s^{{1/ q}-1}\kappa, &\quad
$1\leq\ell\leq K$,}
\nonumber
\\
&& \bT_\varepsilon\dvtx \exists\bigl\{\lambda^k\in
\bbr^m_+,\mu^k\geq0\bigr\}_{k=1}^K\qquad
\forall(k\leq K, t\leq n_k)
\\
&& \qquad\cases{\displaystyle \ErfInv\biggl(\frac{\varepsilon}{2Kn_k} \biggr)
\sum_{t=1}^{n_k}
\bigl\|D^Th^{kt}\bigr\|_2+\frac{1}{2}\sum
_i\lambda^k_i+\frac
{1}{2}
\mu^k\leq\rho,\vspace*{2pt}\cr
\lleft[\matrix{\operatorname{Diag}\bigl\{\lambda^k\bigr
\}&H^T[K]E
\cr
E^TH[k]&\mu^k
I_{n_k}} \rright]\succeq0.}\nonumber
\end{eqnarray}

\subsection{Tractable sufficient conditions: Limits of
performance}\label{seclimitsofperformance}
Consider the situation where all the norms $\|\cdot\|_{(k)}$ are $\|
\cdot\|_r$, with $r\in\{1,2,\infty\}$. A natural question about
verifiable sufficient conditions for a pair $(H, L_\infty(\cdot))$
to satisfy $\bQ_{s,q}(\kappa)$ is, what are the ``limits of
performance'' of these sufficient conditions? Specifically, how large
could be the range of $s$ for which the condition can be satisfied by
at least one contrast matrix? Here is a partial answer to this question:
%
\begin{proposition}\label{proplimits} Let $A$ be an $m\times n$
sensing matrix which is ``essentially nonsquare,'' specifically, such
that $2m\leq n$, let the r.s. $\cS$ be such that $B=I_n$, and let
$n_k=d$, $\|\cdot\|_{(k)}=\|\cdot\|_r$, $1\leq k\leq K$, with $r\in\{
1,2,\infty\}$. Whenever an $m\times n$ matrix $H$ and $n\times n$
matrix $V$ satisfy the conditions
%
\begin{eqnarray}
\label{satiscond}
&I=V+H^TA\quad\mbox{and}&\nonumber\\[-8pt]\\[-8pt]
&\displaystyle \max_{1\leq\ell\leq K}\bigl\|
\bigl[\bigl\| V^{1\ell}\bigr\|_{r,r};\bigl\|V^{2\ell}
\bigr\|_{r,r};\ldots;\bigl\|V^{K\ell}\bigr\|_{r,r}\bigr]
\bigr\|_{s,q} \leq\frac{1}{2}s^{{1/ q}-1}&\nonumber
\end{eqnarray}
[cf. (\ref{equationagain}), (\ref{defnuV}) and (\ref
{goodcase})] with $q\geq1$, one has
%
\begin{equation}
\label{squarerootbound}
s\leq\frac{3\sqrt{m}}{2\sqrt{d}}.
\end{equation}
\end{proposition}

\subsubsection*{Discussion} Let the r.s. $\cS$ in question be the
same as in Proposition~\ref{proplimits}, and let $m\times n$ sensing
matrix $A$ have $2m\leq n$. Proposition~\ref{proplimits} says that in
this case, the verifiable sufficient condition, stated by Proposition
\ref{propnec}, for satisfiability of $\bQ_{s,q}(\kappa)$ with
$\kappa<1/2$ has rather restricted scope---it cannot certify the
satisfiability of $\bQ_{s,q}(\kappa)$, $\kappa\leq1/2$, when $s\geq
\frac{3\sqrt{m}}{2\sqrt{d}}$. Yet, the condition $\bQ_{s,q}(\kappa)$
\textit{may be satisfiable} in a much larger range of values of $s$. For
instance, when the r.s. in question is the standard one and $A$ is a
random Gaussian $m\times n$ matrix, the matrix $A$ satisfies, with
overwhelming probability as $m,n$ grow, the RIP$(\frac{1}{ 5},s)$
condition for $s$ as large as $O(1)m/\sqrt{\ln(n/m)}$ (cf.
\cite{CandesTaorip05}). By Proposition~\ref{appropRIPBlock}, this
implies that $(\frac{5}{4}A, \|\cdot\|_\infty)$ satisfies the condition
$\bQ_{s,2}(\frac{1}{4})$ in essentially the same large range of $s$.
There is, however, an important case where the ``limits of
performance'' of our \textit{verifiable sufficient condition} for the
satisfiability of $\bQ_{s,q}(\kappa)$ implies severe restrictions on
the range of values of $s$ in which the ``true'' condition $\bQ
_{s,q}(\kappa)$ is satisfiable---this is the case when $q=\infty$ and
$r=\infty$. Combining Propositions~\ref{appropTract} and \ref
{proplimits}, we conclude that \textit{in the case of r.s. from
Proposition}~\ref{proplimits} \textit{with $r=\infty$ and}\vadjust{\goodbreak} ``\textit{sufficiently
nonsquare}'' ($2m\leq n$) $m\times n$ \textit{sensing matrix $A$}, \textit{the associated
condition $\bQ_{s,\infty}(\frac{1}{2})$ cannot be satisfied when
$s>\frac{3\sqrt{m}}{2\sqrt{d}}$.}
\subsection{Tractable sufficient conditions and Mutual
Block-Incoherence}\label{secmutual}
We have mentioned in the \hyperref[secintro]{Introduction} that, to the best of our
knowledge, the only previously proposed verifiable sufficient condition
for the validity of $\ell_1$ block recovery is the ``mutual block
incoherence condition'' of \cite{EldarKuppingerBolcskei10}. Our
immediate goal is to show that this condition is covered by Proposition
\ref{propnec}.

Consider an r.s. with $B=I_n$ and with $\ell_2$-norms in the role of
$\|\cdot\|_{(k)}$, $1\leq k\leq K$, and let the sensing matrix $A$ in
question be partitioned as $A=[A[1],\ldots,A[K]]$, where $A[k]$ has $n_k$
columns. Let us define the \textit{mutual block-incoherence} $\mu$ of $A$
w.r.t. the r.s. in question as follows:
%
\begin{equation}\label{MatrixH}\quad
\mu=\mathop{\max_{1\leq k,\ell\leq K,}}_{k\neq\ell} \sigma_{\max} \bigl
(C_k^{-1}A^T[k]A[
\ell] \bigr) \qquad\bigl[\mbox{where }C_k:=A^T[k]A[k]\bigr],
\end{equation}
provided that all matrices $C_k$, $1\leq k\leq K$, are nonsingular,
otherwise $\mu=\infty$. Note that in the case of the standard r.s.,
the just defined quantity is nothing but the standard mutual
incoherence known from the Compressed Sensing literature (see, e.g.,
\cite{DonElTem06}).

In \cite{EldarKuppingerBolcskei10}, the authors consider the same
r.s. and assume that $n_k=d$, $1\leq k\leq K$, and that the columns of
$A$ are of unit $\|\cdot\|_2$-norm. They introduce the quantities
%
\begin{eqnarray}
\label{eqYonina17} \nu&=&\max_{1\le k\leq K} \max_{1\leq j\neq j'\leq K}\bigl|
\Col_j^T\bigl[A[k]\bigr]\Col_{j^\prime
}\bigl[A[k]
\bigr]\bigr|,
\nonumber\\[-8pt]\\[-8pt]
\mu_B&=&\frac{1}{ d}\mathop{\max_{1\leq k,\ell\leq K,}}_{k\neq\ell
}\sigma_{\max}
\bigl(A^T[k]A[\ell]\bigr)
\nonumber
\end{eqnarray}
and prove that an appropriate version of block-$\ell_1$ recovery
allows to recover exactly every $s$-block-sparse signal $x$ from the
noiseless observations $y=Ax$, provided that
%
\begin{equation}
\label{eqYonina55} 1-(d-1)\nu>0 \quad\mbox{and}\quad s<\chi:=\frac{1-(d-1)\nu+
d\mu_B}{2\,d\mu_B}.
\end{equation}

The following observation is almost immediate:
%
\begin{proposition}\label{propYonina} Given an $m\times n$ sensing
matrix $A$ and an r.s. $\cS$ with \mbox{$B=I_n$}, $\|\cdot\|_{(k)}=\|\cdot\|
_2$, $1\leq k\leq K$, let $A=[A[1],\ldots,A[K]]$ be the corresponding
partition of $A$.

\begin{longlist}
\item
Let $\mu$ be the mutual block-incoherence of $A$ w.r.t. $\cS
$. Assuming $\mu<\infty$, we set
%
\begin{eqnarray}
\label{BkAkH}
H=\frac{1}{1+\mu}\bigl[A[1]C_1^{-1},A[2]C_2^{-1},\ldots,A[K]C_K^{-1}
\bigr]
\quad\mbox{where } C_k=A^T[k]A[k].\hspace*{-40pt}
\end{eqnarray}
Then the contrast matrix $H$ along with the matrix $I_n-H^TA$ satisfies
condition (\ref{equationagain}) (where $B=I_n$) and condition
(\ref{sothatthe}) with $q=\infty$ and
\[
\kappa=\frac{\mu s}{1+\mu}.
\]
As a result, applying Proposition~\ref{propnec}, we conclude that whenever
%
\begin{equation}
\label{ssosmall}
s<\frac{1+\mu}{2\mu},
\end{equation}
the pair $(H,L_\infty(\cdot))$ satisfies $\bQ_{s,\infty}(\kappa)$
with $\kappa=\frac{\mu s}{1+\mu}<1/2$.

\item Suppose that $n_k=d, k=1,\ldots,K$, and let the quantities
$\nu$ and $\mu_B$ defined in
(\ref{eqYonina17}) satisfy the relations (\ref{eqYonina55}).
Then\vspace*{1pt} the mutual block-incoherence of $A$ w.r.t. the r.s. in question
does not exceed $\bar{\mu}=\frac{d\mu_B}{1-(d-1)\nu}$.
Furthermore,\vspace*{1pt} we have $\frac{1+\bar{\mu}}{2\bar{\mu}}=\chi$, and
(\ref{ssosmall}) holds, and thus ensures that the contrast $H$,
as defined in (\ref{BkAkH}), and $L_\infty(\cdot)$
satisfy $\bQ_{s,\infty}(\kappa)$ with some $\kappa<\frac{1}{2}$.
\end{longlist}
\end{proposition}

Let $A=[A_{ij}]\in\bbr^{m\times n}$ be a random matrix with i.i.d.
entries $A_{ij}
\sim\break
\cN(0,m^{-1})$. We have the following simple result.
%
\begin{proposition}\label{MuGauss}
Assume that $B=I_n$, $n_k=d$ and $\|\cdot\|_{(k)}=\|\cdot\|_2$ for
all $k$. There are absolute constants $C_1, C_2<\infty$ (the
corresponding bounds are provided in Section \textup{A.10} of the
supplementary article \cite{JKNP-suppl}) such that if $m\ge C_1(d+\ln
(n))$, then the mutual block-incoherence $\mu$ of $A$ satisfies
%
\begin{equation}\label{muboundG}
\mu\le C_2\sqrt{\frac{d+\ln(n)}{ m}}
\end{equation}
with probability at least $1-\frac{1}{ n}$.
\end{proposition}
The bound (\ref{muboundG}), along with Proposition \ref
{propYonina}(i), implies that when $A$ is a Gaussian matrix, all
block-norms are the $\ell_2$-norms and all $n_k=d$ with $d$ ``large
enough'' [such that $d^{-1}\ln n=O(1)$], the verifiable sufficient
condition for $\bQ_{s,\infty}(\frac{1}{3})$ holds with overwhelming
probability for $s=O(\sqrt\frac{m}{ d})$. In other words, in this case
the (verifiable!) condition $Q_{s,\infty}(\kappa)$ attains (up to an
absolute factor) the limit of performance stated in Proposition \ref
{proplimits}.

\section{Matching pursuit algorithm for block recovery}\label{secNEMP}
The Matching Pursuit algorithm for block-sparse recovery is motivated
by the desire to provide a reduced complexity alternative to the
algorithms using $\ell_1$-minimization. Several implementations of
Matching Pursuit for block-sparse recovery have been proposed in the
Compressed Sensing literature \cite
{BaraniukCevher10,Ben-HaimEldar10,EldarKuppingerBolcskei10,EldarMishali09}.
In this section we aim to show that a pair $H,V$ satisfying (\ref
{equationagain}) and (\ref{goodcase}) where $\kappa<1/2$ [and thus,
by Proposition~\ref{propnec}, such that $(H,L_\infty(\cdot))$
satisfies $\bQ_{s,\infty}(\kappa)$] can be used to design a specific
version of the Matching Pursuit algorithm which we refer to as the
\textit{non-Euclidean Block Matching Pursuit} (\textit{NEBMP}) \textit{algorithm} for
block-sparse recovery.\vadjust{\goodbreak}

We fix an r.s. \mbox{$\cS=(B, n_1,\ldots,n_K, \|\cdot\|_{(1)},\ldots,\|\cdot\|
_{(K)})$} and assume
that the block norms $\|\cdot\|_{(k)}$, $k=1,\ldots,K$, are either $\|
\cdot\|_\infty$- or $\|\cdot\|_2$-norms. Furthermore, we suppose
that the matrix $B$ is of full row rank, so that, given $z\in\bbr^N$,
one can compute $x$ such that $z=Bx$ [e.g., $x=B^+z$ where
$B^+=B^T(BB^T)^{-1}$ is the pseudo-inverse of $B$]. Let the noise $\xi
$ in the observation $y=Ax+u+\xi$ be Gaussian, $\xi\sim\cN(0,D)$,
$D\in\bbr^{m\times m}$ is known. Finally, we assume that we are in
the situation of Section~\ref{suffconddisc}, that is, we have at our disposal
an $m\times N$, $N=n_1+\cdots+n_K$, matrix $H$, an $N\times N$ block
matrix $V=[V^{k\ell}
\in\bbr^{n_k\times n_\ell}]_{k,\ell=1}^K$, a $\bar{\gamma}>0$ and
$\rho\geq0$ such that
%
\begin{eqnarray}
\label{conditionNEMP} &&\mbox{\textup{(a)}}\quad B=VB+H^TA,
\nonumber
\\
&&\mbox{\textup{(b)}}\quad\bigl\|V^{k\ell}\bigr\|_{(\ell,k)} =[\Omega]_{k,\ell}\leq\bar{
\gamma}\qquad \forall k,\ell\leq K,
\\
&&\mbox{\textup{(c)}}\quad\Prob_\xi\bigl\{\Xi^+:=\bigl\{\xi\dvtx  L_\infty
\bigl(H^T[u+\xi]\bigr)\leq\rho\ \forall u\in\cU\bigr\}\bigr\}\geq1-
\varepsilon.
\nonumber
\end{eqnarray}

Given observation $y$, a positive integer $s$ and a real $\upsilon\geq
0$ [$\upsilon$ is our guess for an upper bound on $L_1(Bx-[Bx]^s)$],
consider Algorithm~\ref{algogroupMP}.
\begin{algorithm}[t]
\caption{Non-Euclidean Block Matching Pursuit}\label{algogroupMP}
\begin{quote}
1. \textit{Initialization}: Set $v^{(0)}=0$,
$\alpha_0=\frac{L_{s,1}(H^Ty)+s\rho+\upsilon}{1-s\bar{\gamma}}.$

2.
\textit{Step $k$, $k=1,2,\ldots\,$}: Given $v^{(k-1)}\in\bbr^n$ and
$\alpha_{k-1}\geq0$, compute
\begin{enumerate}[2.1.2.\ ]
\item[2.1.] $g=H^T(y-Av^{(k-1)})$ and vector $\Delta=[\Delta
[1],\ldots,\Delta[K]]\in\bbr^N$ by setting for $j=1,\ldots,K$:
%
\begin{eqnarray}
\label{project} \hspace*{-30pt}\Delta[j]&=&\frac{g[j]}{\|g[j]\|_2}\bigl[\bigl\|g[j]\bigr\|_2-\bar{
\gamma}\alpha_{k-1}-\rho\bigr]_+\qquad \mbox{if $\|\cdot\|_{(j)}=\|
\cdot\|_2$};
\nonumber\\[-8pt]\\[-8pt]
\hspace*{-30pt}\Delta_{ji}&=&\sign(g_{ji})\bigl[|g_{ji}|-\bar{
\gamma}\alpha_{k-1}-\rho\bigr]_+,\qquad 1\leq i \leq n_j,\qquad
\mbox{if $
\|\cdot\|_{(j)}=\|\cdot\|_\infty$},\hspace*{-50pt}
\nonumber
\end{eqnarray}
where $w_{ji}$ is $i$th entry in $j$th block of a representation vector
$w$ and $[a]_+=\max\{a,0\}$.

\item[2.2.]
Choose $v^{(k)}$ such that $B(v^{(k)}-v^{(k-1)})=\Delta$, set
%
\begin{equation}
\label{finitedif} \alpha_k= 2s\bar{\gamma}\alpha_{k-1}+2s
\rho+\upsilon.
\end{equation}
%
and loop to step $k+1$.
\end{enumerate}
3. \textit{Output}: The approximate solution found after $k$
iterations is $v^{(k)}$.
\end{quote}
\end{algorithm}
Its convergence analysis is based upon the following:
%
\begin{lemma}\label{PropgroupNEMP}
In the situation of
(\ref{conditionNEMP}), let $s\bar{\gamma}<1$.
Then whenever $\xi\in\Xi^+$, for every $x
\in\bbr^n$ with $L_1(Bx-[Bx]^s)\leq\upsilon$ and every $u \in\cU
$, the following holds true.

When applying Algorithm~\ref{algogroupMP}
to $y=Ax+u+\xi$, the resulting approximations $Bv^{(k)}$ to $Bx$ and
the quantities $\alpha_k$ for all $k$ satisfy the relations
\begin{eqnarray*}
&&\mbox{\textup{(a$_{k}$)}}\quad\mbox{for all $1\leq j\leq K$}\qquad\bigl\|\bigl(Bv^{(k)}-Bx
\bigr)[j]\bigr\|_{(j)}\le\bigl\|(Bx)[j]\bigr\|_{(j)},
\\
&&\mbox{\textup{(b$_{k}$)}}\quad L_1\bigl(Bx-Bv^{(k)}\bigr) \leq
\alpha_{k} \quad\mbox{and}\quad L_\infty\bigl(Bx-Bv^{(k+1)}
\bigr) \leq2\bar{\gamma}\alpha_k+2\rho.
\end{eqnarray*}
\end{lemma}

Note that if $2s\bar{\gamma}<1$, then also $s\bar{\gamma}<1$, so
that Lemma~\ref{PropgroupNEMP} is applicable. Furthermore, in this
case, by (\ref{finitedif}), the sequence $\alpha_k$ converges
exponentially fast to the limit $\alpha_\infty:=\frac{2s\rho
+\upsilon
}{1-2s\bar{\gamma}}$:
\[
L_1\bigl(Bv^{(k)}-Bx\bigr) \leq\alpha_k=(2s
\bar{\gamma})^k[\alpha_0-\alpha_\infty] +
\alpha_\infty.
\]
Along with the second inequality of (b$_{k}$), this implies the bounds
\[
L_\infty\bigl(Bv^{(k)}-Bx\bigr) \leq2\bar{\gamma}
\alpha_{k-1}+2\rho\leq\frac{\alpha_k }{ s},
\]
and since $L_p(w) \leq L_1(w)^{{1}/{ p}} L_\infty(w)^{({p-1
})/{ p}}$
for $1\leq p \leq\infty$, we have
\[
L_p\bigl(Bv^{(k)}-Bx\bigr) \leq s^{({1-p })/{ p}} \bigl[(2s
\bar{\gamma})^k[\alpha_0-\alpha_\infty] +
\alpha_\infty\bigr].
\]

The bottom line here is as follows.
%
\begin{proposition}\label{CorgroupNEMP}
Suppose that a collection $(H, L_\infty(\cdot),\rho,\bar{\gamma
},\varepsilon)$ satisfies (\ref{conditionNEMP}), and let the
parameter $s$ of Algorithm~\ref{algogroupMP} satisfy $2\kappa:=2s\bar{\gamma}<1$.
Then for all $\xi\in\Xi^+$, $u\in\cU$, $x\in\bbr^n$ such that
$L_1(Bx-[Bx]^s)\leq\upsilon$, Algorithm~\ref{algogroupMP}
as applied to $y=Ax+u+\xi$ ensures that for every $t=1,2,\ldots$ one has
\begin{eqnarray*}
&&L_p\bigl(Bv^{(t)}-Bx\bigr)\\
&&\qquad\leq s^{{1}/{ p}} \biggl[
\frac{2\rho+s^{-1}\upsilon}{1-2\kappa}
+(2\kappa)^{t} \biggl(\frac{ s^{-1} (L_{s,1}(H^Ty) +\upsilon) +\rho}{
1-\kappa} -\frac{2\rho+s^{-1}\upsilon}{1-2\kappa} \biggr)
\biggr] 
\end{eqnarray*}
for all $1\leq p\leq\infty$'s \textup{[cf. (\ref{apeq60})]}.
\end{proposition}
Note that Proposition~\ref{CorgroupNEMP} combined with Proposition
\ref{propYonina} essentially covers the results of \cite
{EldarKuppingerBolcskei10} on the properties of the Matching Pursuit
algorithm for the block-sparse recovery proposed in this reference.

\section{Numerical illustration}\label{secnumerics}
In the theoretical part of this paper we considered the situation where
the sensing matrix $A$ and the r.s. $\cS=(B,n_1,\ldots,n_K,\break \|\cdot
\|_{(1)},\ldots,\|\cdot\|_{(K)})$ were given, and we were interested in
understanding:
\begin{longlist}[(B)]
\item[(A)] whether $\ell_1$ recovery allows to recover the
representations $Bx$ of all $s$-block-sparse signals with a given $s$
in the absence of observation noise, and
\item[(B)] how to choose the best (resulting in the smallest
possible error bounds) pair $(H,\|\cdot\|)$.\footnote{Needless to
say, the results presented so far do not pretend to provide full
answers to these questions. Our verifiable sufficient conditions for
the validity of $\ell_1$ block recovery supply only \textit{lower bounds}
on the largest $s=s_*$ for which the answer to (A) is positive.
Similarly, aside of the case $q=\infty$, $\|\cdot\|_{(k)}=\|\cdot\|
_\infty$, $1\leq k\leq K$, our conditions for the validity of
block-$\ell_1$ recovery are only sufficient, meaning that optimizing
the error bound over $(H,\|\cdot\|)$ allowed by these conditions may
only yield \textit{sub}optimal recovery routines.}
\end{longlist}
Note that our problem setup involves a number of components. While in
typical applications sensing matrix $A$, representation matrix $B$ and
the dimensions $n_1,\ldots,n_K$ of the block vectors may be thought as
given by the ``problem's physics,'' it is not the case for the block
norms $\|\cdot\|_{(k)}$. Their choice (which does affect the $\ell_1$ recovery routines) appears to be unrelated to the model of the data.

\textit{The first goal} of our experiments is to understand how to choose
the block norms in order to validate $\ell_1$ recovery for the largest
possible value of the sparsity parameter $s$; here ``to validate''
means to provide guarantees of small recovery error for all
$s$-block-sparse signals when the observation error is small (which
implies, of course, the exactness of the recovery in the case of
noiseless observation). Here we restrict ourselves to the case of $\ell
_r$-r.s. with $r\in\{1,2,\infty\}$.
By reasons explained in the discussion in Section~\ref{secrecovery},
we consider here only the case of the penalized $\ell_1$ recovery with
$m\times N$ contrast matrix $H$ (where, as always, $N=n_1+\cdots+n_K$), $\|
\cdot\|=L_\infty(\cdot)$,\footnote{These are exactly the pairs
$(H,\|\cdot\|)$ covered by the sufficient conditions for the validity
of $\ell_1$ recovery; see Proposition~\ref{propnec}.} and with
$\lambda=2s$ [see (\ref{apeq4pen})]. Besides this, we assume, mainly
for the sake of notational convenience, that $B=I_n$.

Let us fix $A\in\bbr^{m\times n},B=I_n,K,n_1,\ldots,n_K$
($n_1+\cdots+n_K=n=:N$). By Proposition~\ref{propnec}, for every matrix
$H\in\bbr^{m\times n}$
setting
%
\begin{eqnarray}
\label{numeq2} V&\equiv&\bigl[V^{k\ell}\in\bbr^{n_k\times n_\ell}
\bigr]_{k\ell=1}^K=I-H^TA,\nonumber\\
\Omega^r(H)&=&
\bigl[\bigl\|V^{k\ell}\bigr\|_{r,r}\bigr]_{k,\ell=1}^K,
\nonumber\\[-8pt]\\[-8pt]
\qquad\kappa_1^{r,s}(H)&=&\max_{1\leq\ell\leq K}\bigl\|\Col_\ell
\bigl[\Omega^r(H)\bigr]\bigr\|_{s,1},\nonumber\\
\kappa_\infty^{r,s}(H)&=&s
\max_{1\leq k,\ell\leq K}\bigl[\Omega^r(H)\bigr]_{k,\ell},
\nonumber
\end{eqnarray}
the pair $(H,L_\infty(\cdot))$ satisfies the conditions $\bQ
_{s,q}(\kappa^{r,s}_q(H))$, $q=1$ and $q=\infty$, provided that the
block norms are the $\ell_r$-ones. In particular, when $\kappa
_1^{r,s}(H)<1/2$, the penalized $\ell_1/\ell_r$ recovery [i.e., the
recovery (\ref{apeq4}) with all block norms being the $\ell_r$-ones]
``is valid'' on $s$-block-sparse signals, meaning exactly that this
recovery ensures the validity of the error bounds (\ref
{apeq2626standard}) with $q=\infty$, $\varkappa=\kappa^{r,s}_1$,
$\kappa=\kappa^{r,s}_\infty$ (and, in particular, recovers exactly
all $s$-block-sparse signals when there is no observation noise).

\textit{Our strategy} is as follows. For each value of $r\in\{1,2,\infty
\}$, we consider the convex optimization problem
\[
\min_{H\in\bbr^{m\times n}} \Bigl\{\kappa_1^{r,s}(H):=
\max_{\ell\leq K}\bigl\|\Col_\ell\bigl[\Omega^r(H)\bigr]
\bigr\|_{s,1} \Bigr\},%
\]
find the largest $s=s(r)$ for which the optimal value in this problem
is $<1/2$, and denote by $H^{(r)}, r\in\{1,2,\infty\}$ the
corresponding optimal solution. In addition to these ``marked''
contrast matrices, we consider two more contrasts, $H^{(\mathrm{MI})}$
and $H^{(\mathrm{MBI})}$, based on the mutual
block-incoherence condition and given by the calculation
(\ref{MatrixH}) for the cases of the ``standard'' (1-element blocks in
$x=Bx$) and the actual block structures, respectively.\looseness=1

Now, given the set $\cH=\{H^{(\mathrm{MI})}, H^{(\mathrm{MBI})}, H^{(1)}, H^{(2)}, H^{(\infty)}\}$ of
$m\times n$ candidate contrast matrices, we can choose the ``most
powerful'' penalized $\ell_1/\ell_r$ recovery suggested by $\cH$ as
follows: for every $H\in\cH$ and for every $p\in\{1,2,\infty\}$, we
find the largest $s=s(H,p)$ for which $\kappa_1^{r,p}(H)<1/2$, and
then define the quantity $s_*=s_*(\cH)=\max\{s(H,p)\dvtx  H\in\cH,
p\in\{1,2,\infty\}\}$ along with $H_*\in\cH$ and $p_*\in\{
1,2,\infty\}$ such that $s_*=s(H_*,p_*)$. The penalized $\ell_1/\ell
_{p_*}$ recovery utilizing the contrast matrix $H_*$ and the norm
$L_\infty(\cdot)$ associated with block norms $\|\cdot\|_{p_*}$ of
the blocks is definitely valid for $s=s_*(\cH)$, and this is the
largest sparsity range, as certified by our sufficient conditions, for
the validity of $\ell_1/\ell_r$ recovery, which we can get with
contrast matrices from $\cH$. Note that $s_*\geq\max
[s(1),s(2),s(\infty)]$, that is, the resulting range of values of $s$
is also the largest we can certify using our sufficient conditions,
with \textit{no} restriction on the contrast matrices.

\subsection*{Implementation} We have tested the outlined strategy in
the following problem setup:
\begin{itemize}
\item the sensing matrices $A$ are of \textit{size} $(m=96)\times
(n=128)$, $B=I$ with $K=32$ four-element blocks in $Bx=x$;
\item the \mbox{$96\times128$} \textit{sensing matrices} $A$ are built as
follows: we first draw a matrix at random from one of the following
distributions:

%
\begin{table}
\def\arraystrech{0.9}
\caption{Certified sparsity levels for penalized $\ell
_1/\ell_r$-recoveries for candidate contrast matrices.
For each candidate and each value of $r$ we present in the
corresponding cells\break the triple $s(H,r) |\kappa_1^{r,s(H,r)}(H)
|\kappa_\infty^{r,s(H,r)}(H)$. $\bar{s}(r)$: a computed upper bound
on\break $r$-goodness $s^*(A,r)$ of $A$. Italic: the best sparsity
$s_*(\cH)$ certified by our sufficient conditions for the validity of
penalized recovery}\label{table1}\vspace*{-3pt}
\begin{tabular*}{\tablewidth}{@{\extracolsep{4in minus 4in}}lccd{1.4}cc
ccccc@{}}
\hline
$\bolds{A}$&$\bolds{r}$&\multicolumn{3}{c}{$\bolds{H^{(\mathrm{MI})}}$}&
\multicolumn{3}{c}{$\bolds{H^{(\mathrm{MBI})}}$}
&\multicolumn{3}{c@{}}{$\bolds{H^{(1)}}$}\\
\hline
\texttt{H}&$1$&2&0.4727&0.509&2&0.444&0.460&\textit{3}&0.429&0.429\\
&$2$&2&0.436&0.436&2&0.429&0.429&\textit{3}&0.429&0.429\\
&$\infty$&2&0.473&0.509&2&0.444&0.460&\textit{3}&0429&0.429\\
[3pt]
\texttt{G}&$1$&0&0.000&0.000&0&0.000&0.000&\textit{3}&0.467&0.900\\
&$2$&0&0.000&0.000&1&0.368&0.368&1&0.300&0.300\\
&$\infty$&0&0.000&0.000&0&0.000&0.000&0&0.000&0.000\\
[3pt]
\texttt{R}&$1$&0&0.0000&0.000&0&0.000&0.000&\textit{3}&0.477&0.853\\
&$2$&0&0.000&0.000&1&0.354&0.354&1&0.284&0.284\\
&$\infty$&0&0.000&0.000&0&0.000&0.000&1&0.482&0.482\\
[3pt]
\texttt{T}&$1$&1&0.384&0.384&1&0.399&0.399&\textit{2}&0.383&0.383\\
&$2$&1&0.384&0.384&1&0.399&0.399&\textit{2}&0.383&0.383\\
&$\infty$&1&0.384&0.384&1&0.399&0.399&\textit{2}&0.383&0.383\\
\hline
\end{tabular*}

\begin{tabular*}{\tablewidth}{@{\extracolsep{\fill}}lcccccccc@{}}
$\bolds{A}$&$\bolds{r}$&\multicolumn{3}{c}{$\bolds{H^{(2)}}$}
&\multicolumn{3}{c}{$\bolds{H^{(\infty)}}$}&$\bolds{\bar{s}(r)}$\\
\hline
\texttt{H}&$1$&2&0.487&0.519&\textit{3}&0.429&0.429&4\\
&$2$&\textit{3}&0.429&0.429&\textit{3}&0.429&0.429&3\\
&$\infty$&2&0.487&0.519&\textit{3}&0.429&0.429&3\\
[3pt]
\texttt{G}&$1$&1&0.301&0.301&1&0.489&0.489&5\\
&$2$&\textit{3}&0.447&0.458&2&0.479&0.549&5\\
&$\infty$&1&0.305&0.305&\textit{3}&0.483&0.823&4\\
[3pt]
\texttt{R}&$1$&1&0.291&0.291&1&0.498&0.498&5\\
&$2$&\textit{3}&0.438&0.440&1&0.264&0.264&5\\
&$\infty$&1&0.286&0.286&\textit{3}&0.489&0.739&5\\
[3pt]
\texttt{T}&$1$&\textit{2}&0.383&0.383&\textit{2}&0.383&0.383&3\\
&$2$&\textit{2}&0.383&0.383&\textit{2}&0.383&0.383&3\\
&$\infty$&\textit{2}&0.383&0.383&\textit{2}&0.383&0.383&3\\
\hline
\end{tabular*}   \vspace*{-3pt}
\end{table}

\begin{itemize}
\item\textit{type} \texttt{H}: randomly selected $96\times128$ submatrix of
the $128\times128$ Hadamard matrix,\footnote{The Hadamard matrices
$H_k$ of order $2^k\times2^k$, $k=0,1,\ldots\,$, are given by the
recurrence $H_0=1$, $H_{k+1}=[H_k,H_k;H_k,-H_k]$. They are symmetric
matrices with $\pm1$ entries and rows orthogonal to each other.}
\item\textit{type} \texttt{G}: $96\times128$ matrix with independent $\cN
(0,1)$ entries,
\item\textit{type} \texttt{R}: $96\times128$ matrix with independent
entries taking values $\pm1$ with equal probabilities,
\item\textit{type} \texttt{T}: random $96\times128$ matrix of the structure
arising in \textit{Multi-Task Learning} (see, e.g., \cite{MTL} and
references therein): the consecutive 4-column parts of the matrix are
block-diagonal with four $24\times1$ diagonal blocks with independent
$\cN(0,1)$ entries,
\end{itemize}
and then scale the columns of the selected matrix to have their $\|
\cdot\|_2$-norms equal to 1.
\end{itemize}

\textit{The results} we report describe 4 experiments differing from
each other by the type of the (randomly selected) matrix
$A$.\footnote{As far as our experience shows, the results remain nearly
the same across instances of $A$ drawn from the same distribution, so
that only one experiment for each type of distribution in question
appears to be representative enough.}

\begin{figure}

\includegraphics{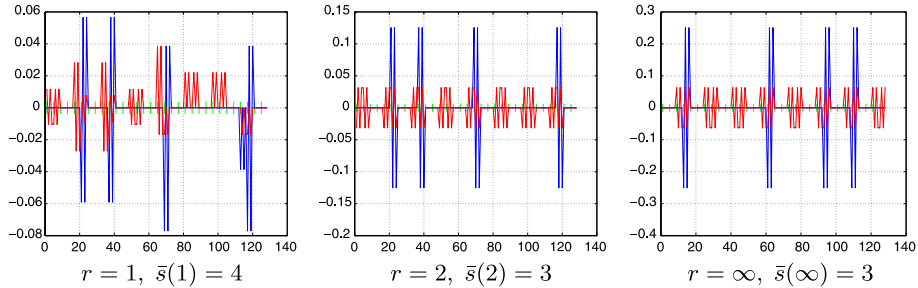}

\caption{``Bad'' $(\bar{s}(r)+1)$-block-sparse signals (blue) and
their $\ell_1/\ell_r$ recoveries (red) from noiseless observations,
\texttt{H}-matrix $A$.}
\label{fig1}
\end{figure}

In Table~\ref{table1}, we display the certified sparsity levels of
penalized $\ell_1/\ell_r$ recoveries for the candidate contrast matrices.
In addition, we present valid \textit{upper bounds} $\bar{s}(r)$ on the
``$r$-goodness''
$s^*(A,r)$ of $A$, defined as the largest $s$ such\vadjust{\goodbreak} that the $\ell_1/\ell
_r$ recovery in the noiseless case recovers \textit{exactly} the
representations of \textit{all} $s$-block-sparse vectors, that is,
\begin{eqnarray*}
s^*(A,r)&=&\max\Biggl\{s\dvtx  x=\mathop{\operatorname{Arg}\min}_{z\in\bbr^n}
\Biggl
\{\sum_{k=1}^K\bigl\|[z]_k
\bigr\|_r\dvtx Az=Ax \Biggr\}
\\
&&\hspace*{93.5pt}\mbox{for all $s$-block-sparse $x$.}\Biggr\}
\end{eqnarray*}
We present on Figure~\ref{fig1} examples of ``bad'' signals [i.e.,
$(\bar{s}(r)+1)$-block-sparse signals which are \textit{not} recovered
correctly by the latter procedure].\footnote{It is immediately seen
that whenever $B$ is of full row rank, the \textit{nullspace property}
``$L_{s,1}(Bx)<\frac{1}{2}L_1(Bx)$ for all $x\in\Ker A$ with $Bx\neq
0$'' is necessary for $s$ to be $\leq s^*(A,\cdot)$. As a result, for
$B$'s of full row rank, $s^*(A,r)$ can be upper-bounded in a manner
completely similar to the case of the standard r.s.; see \cite{JNCS}, Section
4.1.}

On the basis of this experiment we can make two tentative conclusions:
\begin{itemize}
\item the $\ell_1/\ell_2$ recovery with the contrast matrix $H^{(2)}$
and the $\ell_1/\ell_\infty$ recovery with the contrast matrix
$H^{(\infty)}$
were able to certify the best levels of allowed sparsity (when compared
to other candidate matrices from $\cH$);
\item in our experiments, the upper bounds $\bar{s}(r)$ on the
$r$-goodness $s^*(A,r)$ of $A$ are close to the corresponding certified
lower bounds $s_*(\cH,r)=\max_{H\in\cH} s(H,r)$.
\end{itemize}
%

\textit{Numerical evaluation of recovery errors.} The objective of the
next experiment is to evaluate the accuracy of penalized $\ell_1/\ell
_r$ recoveries in the noisy setting. As above, we consider the contrast
matrices from $\cH=\{H^{(\mathrm{MI})}, H^{(\mathrm{MBI})},\break
H^{(1)},H^{(2)}, H^{(\infty)}\}$. Note that it is possible to
improve the error bound by optimizing it over $H$ as it was done in
Section~\ref{suffconddisc}. In the experiments to be reported this
additional optimization, however, did not yield a significant
improvement (which perhaps reflects the ``nice conditioning'' of the
sensing matrices we dealt with), and we do not present the simulation
results for optimized contrasts here:

\begin{itemize}
\item We ran four series of simulations corresponding to the four
instances of the sensing matrix $A$ we used. The series associated with
a particular $A$ was as follows:
\item Given $A$, we associate with it the five aforementioned
candidate contrast matrices from $\cH$. Combining these matrices with
3 values of $r$ ($r=1,2,\infty$), we get 15 recovery routines. In
addition to these 15 routines, we also included the block Lasso
recovery as described in \cite{Lounicietal10}. In our notation,
this recovery is (cf. \cite{Lounicietal10}, (2.2))
\[
\widehat{x}_{\mathrm{Lasso}}(y)\in\Argmin_z \Biggl\{
\frac{1}{ m}\|Az-y\|_2^2+2\sum
_{k=1}^K\lambda_k\bigl\|z[k]\bigr\|_2
\Biggr\}
\]
($z[k]$, $1\leq k\leq K$, are the blocks in $z=Bz$), with the penalty
coefficients $\lambda_k$ chosen according to the equality version of
the relations in \cite{Lounicietal10}, Theorem~3.1, used with $q=2$.
\end{itemize}
Each of the 16 resulting recovery routines was tested on two samples,
each containing 100 randomly generated recovery problem instances. In
each problem instance the true signal was randomly generated with $s$
nonzero blocks, and the observations were corrupted by pure Gaussian
white noise: $y=Ax+\sigma\xi$, $\xi\sim\cN(0,I)$. In the first
sample, $s$ was set to the best value $s_*(\cH)$ of block sparsity we
were able to certify; in the second, $s=2s_*(\cH)$ was used. The
parameter $\lambda$ of the penalized recoveries was set to $2s$ (and
thus was tuned to the actual sparsity of test signals). In both
samples, we used $\sigma=0.001$.

\begin{table}
\tablewidth=240pt
\caption{Ratings of recovery routines}\label{table2a}
\begin{tabular*}{\tablewidth}{@{\extracolsep{\fill}}lcccccc@{}}
\hline
$\bolds{r}$&$\bolds{H^{(\mathrm{MI})}}$&$\bolds{H^{(\mathrm{MBI})}}$
&$\bolds{H^{(1)}}$&$\bolds{H^{(2)}}$&$\bolds{H^{(\infty)}}$&\textbf{Lasso}\\
\hline
1&0.30&0.20&0.53&0.60&0.54&N/A\\
2&0.76&0.51&0.75&0.79&0.75&0.19\\
$\infty$&0.25&0.18&0.44&0.48&0.44&N/A\\
\hline
\end{tabular*}     \vspace*{-3pt}
\end{table}

We compare the recovery routines on the basis of their \textit{ratings}
computed as follows: given a recovery problem instance from the sample,
we applied to it every one of our 16 recovery routines and measured the
16 resulting $\|\cdot\|_\infty$-errors. Dividing the smallest of
these errors by the error of a given routine, we obtain ``the rating''
of the routine in this particular simulation. Thus, all ratings are
$\leq1$; and the routine which attains the best $\|\cdot\|_\infty
$ recovery error for the current data is rated ``1.0.'' For the
remaining routines, the closer to 1 is the rating of the routine, the
closer is the routine to the ``winner'' of the current simulation. The
final rating of a given recovery routine is its average rating over all
$800=4\times2\times100$ recovery problem instances processed in the
experiment.

\begin{figure}

\includegraphics{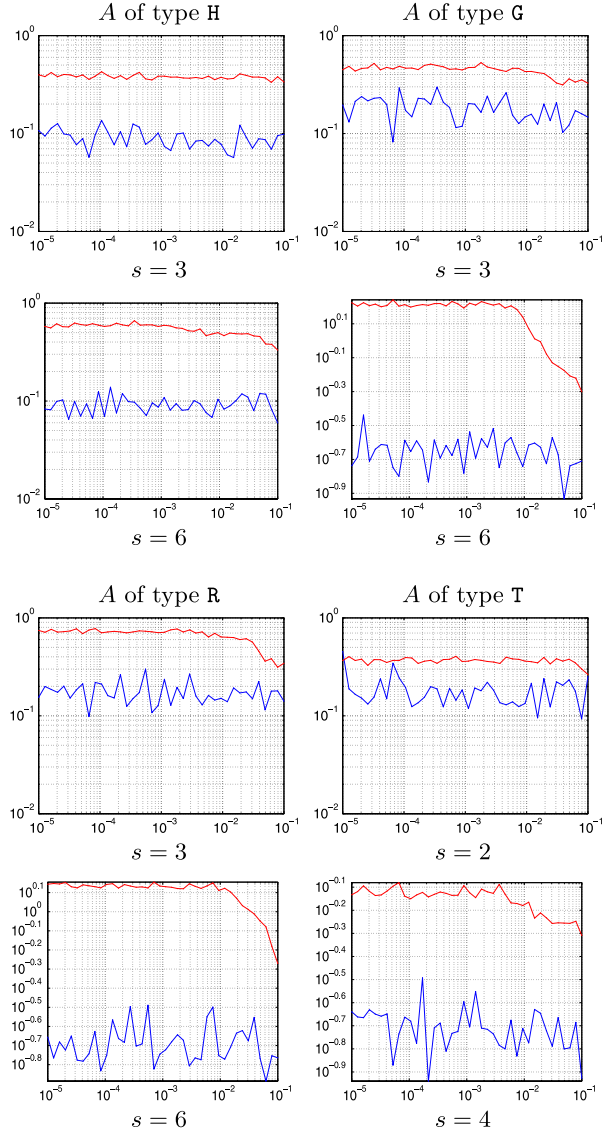}

\caption{Average over 40 experiments ratio of $\|\cdot\|_\infty
$ recovery error to $\sigma$ vs. $\sigma$. In blue: $\ell_1/\ell_2$ recovery with $H=H^{(2)}$; in red: Lasso recovery.}
\label{fig2}
\end{figure}

The resulting ratings are presented in Table~\ref{table2a}. The
``winner'' is the routine associated with $r=2$ and $H=H^{(2)}$.
Surprisingly, the second best routine is associated with the same $r=2$
and the simplest contrast $H^{(\mathrm{MI})}$, an
outsider in terms of the data presented in Table~\ref{table1}. This
inconsistency may be explained by the fact that the data in Table \ref
{table1} describe the guaranteed worst-case behavior of our recovery
routines, which may be quite different from their ``average behavior,''
reflected by Table~\ref{table2a}. Our tentative conclusion on the
basis of the data from Tables~\ref{table1} and~\ref{table2a} is that
the penalized $\ell_1/\ell_2$ recovery associated with the contrast
matrix $H^{(2)}$ may be favorable when recovery guarantees are to be
associated with good numerical performance.

The above comparison was carried out for $\sigma$ set to 0.001. The
conducted experiments show that for the routines in question and our
purely Gaussian model of observation errors, the recovery errors are,
typically, proportional to $\sigma$. This is illustrated by the plots
on Figure~\ref{fig2} where we traced the average (over 40 experiments
for every grid value of $\sigma$) \textit{signal-to-noise ratio} (the
ratio of the $\|\cdot\|_\infty$-error of the recovery to $\sigma$)
of our favorable recovery ($r=2$, $H=H^{(2)}$) and the corresponding
performance figure for block Lasso.

\begin{supplement}
\stitle{Supplement to ``Accuracy guaranties for $\ell_1$ recovery of block-sparse signals''}
\slink[doi]{10.1214/12-AOS1057SUPP} 
\sdatatype{.pdf}
\sfilename{aos1057\_supp.pdf}
\sdescription{The proofs of the results stated in the paper and the
derivations for Section~\ref{suffconddisc} are provided in the
supplementary article \cite{JKNP-suppl}.}
\end{supplement}


\printaddresses

\end{document}